\documentclass[12pt]{amsart}

\calclayout

\usepackage[utf8]{inputenc}
\usepackage{amsmath,amsthm,amssymb}
\usepackage{colonequals}
\usepackage{enumitem}

\usepackage[colorlinks=true,allcolors=blue,pdftitle={Metastability of the Meissner State}]{hyperref}
\usepackage[initials,nobysame,alphabetic,msc-links,backrefs]{amsrefs}
\usepackage[foot]{amsaddr}

\newtheorem{theorem}{Theorem}[section]
\newtheorem{lemma}{Lemma}[section]
\newtheorem{lemmaA}{Lemma}

\newtheorem{prop}{Proposition}[section]
\newtheorem{propA}{Proposition}

\newtheorem{remark}{Remark}[section]
\newtheorem{definition}{Definition}[section]

\numberwithin{equation}{section}

\DeclareMathOperator{\diver}{div}
\DeclareMathOperator{\curl}{curl}

\newcommand{\R}{\ensuremath{\mathbb{R}}}
\newcommand{\C}{\ensuremath{\mathbb{C}}}
\newcommand{\HH}{\ensuremath{\mathcal{H}}}
\newcommand{\GG}{\ensuremath{\mathfrak{G}}}
\newcommand{\RR}{\ensuremath{\mathfrak{R}}}

\newcommand{\ep}{\varepsilon}
\newcommand{\de}{\delta}
\newcommand{\ex}{\mathrm{ex}}
\newcommand{\loc}{\mathrm{loc}}
\newcommand{\lep}{{|\log\ep|}}

\renewcommand{\u}{{\mathsf{u}}}
\newcommand{\A}{{\mathsf{A}}}
\renewcommand{\H}{{\mathsf{H}}}

\newcommand{\pmin}{\ensuremath{\rho_\ep}}
\newcommand{\Bfield}{\ensuremath{B_\ep^0}}

\newcommand{\fenthreedOmega}{\ensuremath{F_{\ep,\Omega}}}
\newcommand{\fenweightthreed}{\ensuremath{F_{\ep,\pmin}}}
\newcommand{\glthreed}{\ensuremath{GL_\ep}}

\newcommand{\phase}{\ensuremath{e^{ih_\ex\phi_\ep^0}}}
\newcommand{\meisconf}{\ensuremath{\pmin\phase,h_\ex A_\ep^0}}
\newcommand{\meisdecom}{\ensuremath{\pmin u\phase,A+h_\ex A_\ep^0}}

\newcommand{\ACouGauge}{\ensuremath{\A^C}}
\newcommand{\AUniGauge}{\ensuremath{\widetilde{\A}}}
\newcommand{\uCouGauge}{\ensuremath{\u^C}}
\newcommand{\uUniGauge}{\ensuremath{\widetilde{\u}}}

\newcommand{\pzero}{{p_0}}

\newcommand{\minsp}{H^1(\Omega,\C)\times[A_{\ex}+H_{\curl}]}

\title[Metastability of the Meissner state with pinning]{Metastability of the Meissner state in the three-dimensional Ginzburg--Landau model with pinning}

\author[Carlos Rom\'{a}n]{Carlos Rom\'{a}n}
\address{Facultad de Matem\'aticas \& Instituto de Ingenier\'ia Matem\'atica y Computacional, Pontificia Universidad Cat\'olica de Chile, Vicu\~na Mackenna 4860, 7820436 Macul, Santiago, Chile}
\email{carlos.roman@uc.cl}

\date{October 30, 2026}



\thanks{Funding information: ANID FONDECYT 1231593.}

\begin{document}
	
\begin{abstract}
	We study the persistence of the Meissner state in the three-dimensional magnetic Ginzburg--Landau model with a spatially dependent pinning term. Factoring out the equilibrium profile associated with the pinning term reduces the problem to the analysis of a weighted Ginzburg--Landau functional. In the extreme type-II regime, we construct locally minimizing vortexless configurations for applied fields extending far beyond the first critical field. We obtain quantitative estimates showing that these configurations remain close, up to gauge transformations, to a Meissner-type state determined by a weighted London problem. We also prove the uniqueness, up to gauge transformations, of the locally minimizing Meissner-type branch within a quantitatively controlled class of local minimizers. A central difficulty is to derive elliptic estimates that remain uniform in the Ginzburg--Landau parameter in the presence of the resulting spatially dependent weight. Under a small-BMO assumption on this weight, we obtain the required regularity and extend to the three-dimensional pinned setting the metastability picture previously established for the homogeneous model.
\end{abstract}
\maketitle 

%%%%%%%%%%%%%%%%%%%%%%%%%%%%%%%%%%%%%%%%%%%%%%%%%
\section{Introduction}
%%%%%%%%%%%%%%%%%%%%%%%%%%%%%%%%%%%%%%%%%%%%%%%%%

The magnetic Ginzburg--Landau model \cite{GinLan} provides a phenomenological description of superconducting materials under applied magnetic fields. Its configurations depend on a complex-valued order parameter, which describes the local superconducting state, and on a magnetic vector potential for the induced magnetic field. A fundamental parameter of the model is the Ginzburg--Landau parameter $\kappa$, which distinguishes between type-I and type-II superconductors. Throughout this work, we write $\ep=\kappa^{-1}$ and consider the regime $\ep\to0$, corresponding to extreme type-II superconductors.

At low applied fields, a superconducting sample is expected to be in the Meissner state, characterized by the expulsion of the magnetic field from the interior of the material. As the intensity $h_{\ex}$ of the applied field increases, this state eventually ceases to be globally minimizing and configurations containing vortices become energetically favorable. Vortices are defects in the superconducting state that carry quantized magnetic flux; in three dimensions, they concentrate along curves and are therefore referred to as vortex lines or filaments. The threshold at which they first become energetically favorable is known as the first critical field and is denoted by $H_{c_1}$. In the extreme type-II regime, its leading order is proportional to $|\log\ep|$.

The first critical field is, however, an energetic threshold and does not necessarily coincide with the loss of stability of the Meissner state. Although vortex configurations have lower energy above $H_{c_1}$, a locally minimizing vortexless solution may persist for substantially larger applied fields. The field at which this branch loses stability is known as the \emph{superheating field}, denoted by $H_{\mathrm{sh}}$.

This distinction was made rigorous in the three-dimensional homogeneous setting in \cite{Rom-CMP}. More precisely, it was proved that below the first critical field the Meissner solution is, up to gauge transformations, the unique global minimizer of the Ginzburg--Landau energy. Beyond this threshold, a Meissner-type solution persists as a local minimizer at least for $h_{\ex}=o(\ep^{-\frac13})$, and is unique, up to gauge transformations, at least for $h_{\ex}=o(\ep^{-\frac14})$. Thus, when the applied field is gradually increased, vortex nucleation is expected to be delayed beyond $H_{c_1}$ and to occur only when the vortexless branch becomes unstable. This loss of stability is associated with the superheating field, which is expected to be of order $\ep^{-1}$; see, for instance, \cite{chapman}.

The purpose of the present work is to investigate to what extent this picture remains valid for the three-dimensional Ginzburg--Landau functional in the presence of material inhomogeneities described by a spatially dependent pinning term.

\subsection{The Ginzburg--Landau model with pinning}

The performance of type-II superconductors is strongly affected by the motion of vortices, which produces dissipation and limits the current that the material can carry while remaining superconducting. Material inhomogeneities may restrict this motion by acting as pinning sites for vortices. This provides one of the main motivations for introducing spatially dependent coefficients into the Ginzburg--Landau model. Such an approach goes back to Likharev \cite{Likha}, who proposed to account for spatial variations of the superconducting properties by allowing the coefficients in the energy functional to depend on the position.

Vortex pinning has been the subject of considerable interest in the physics of superconductivity. The review \cite{Kwok2016}, for instance, discusses how different types of defects and inclusions give rise to complex pinning landscapes, which influence vortex configurations and play an important role in enhancing the critical current, namely the maximal current that can be carried without losing the superconducting behavior. In a more microscopic direction, experimental work on FeSe-based superconductors has provided direct evidence that local defects can act as pinning centers for vortices \cite{phys-x-2024}. These works provide physical motivation for considering spatially dependent pinning terms in the Ginzburg--Landau functional rather than restricting the analysis to the classical homogeneous model. For a more detailed discussion of the physics underlying vortex pinning and its implications, see, for instance, \cites{Likha,Tin,chapman-richardson-pinning} and the references therein.

\medskip

After nondimensionalization of the physical constants, the Ginzburg--Landau functional with pinning, which models an inhomogeneous superconducting sample in an applied magnetic field at a fixed temperature below the critical one, is given by
\begin{equation}\label{GLenergy}
	\glthreed(\u,\A)=\frac12\int_\Omega \left(|\nabla_{\A} \u|^2+\frac{1}{2\ep^2}(a_\ep(x)-|\u|^2)^2\right)+\frac12\int_{\R^3}|\H-H_\ex|^2.
\end{equation}
Here, $\Omega$ is a bounded domain of $\R^3$, which we assume to be simply connected with $C^2$ boundary. The function $\u:\Omega\rightarrow\mathbb{C}$ is the \emph{order parameter}; its squared modulus, which represents the relative density of Cooper pairs of superconducting electrons in the microscopic BCS theory \cite{BCS}, indicates the local state of the superconductor. The vector field $\A:\R^3\rightarrow\R^3$ is the electromagnetic vector potential of the induced magnetic field $\H=\curl\A$, and $\nabla_\A$ denotes the covariant gradient $\nabla-i\A$. The external, or applied, magnetic field is assumed to be of the form $H_\ex=h_{\ex} H_{0,\ex}$, where $H_{0,\ex}\in L^\infty(\R^3,\R^3)$ is a fixed divergence-free vector field and $h_{\ex}\geq0$ is a tunable parameter representing the intensity of the external field.

Finally, $a_\ep$ is a function that accounts for inhomogeneities in the material. We assume that $a_\ep\in L^\infty(\Omega)$ and that it takes values in $[b,1]$, where $b\in(0,1)$ is a constant independent of $\ep$. Regions where $a_\ep=1$ correspond to sites without inhomogeneities, whereas regions where $a_\ep<1$ may act as pinning sites.

\smallskip

The Ginzburg--Landau model is a $\mathbb{U}(1)$-gauge theory, in which all the physically relevant quantities are invariant under the gauge transformations
$$
\u\to\u e^{i\Phi}, \qquad \A\to\A+\nabla\Phi,
$$
where $\Phi$ is any sufficiently regular real-valued function. Two important gauge-invariant quantities are the supercurrent and the vorticity, respectively defined, for any sufficiently regular configuration $(\u,\A)$, by
\begin{equation*}
	j(\u,\A)=(i\u,\nabla_\A\u), \qquad \mu(\u,\A)=\curl j(\u,\A)+\curl\A,
\end{equation*}
where $(\cdot,\cdot)$ denotes the scalar product in $\mathbb{C}$ identified with $\mathbb{R}^2$.

Another gauge-invariant quantity that will play an important role throughout this paper is the free energy restricted to $\Omega$, defined by
$$
\fenthreedOmega(u,A)=\frac12\int_\Omega \left( |\nabla_A u|^2+\frac{1}{2\ep^2}(1-|u|^2)^2+|\curl A|^2\right).
$$

\subsection{A brief historical overview}

The homogeneous Ginzburg--Landau model has a long history of mathematical analysis. The two-dimensional theory was developed in a series of important works; we refer, in particular, to the foundational monographs \cites{BetBreHel,SanSerBook} and the references therein for comprehensive accounts of the subject.

In the two-dimensional homogeneous setting, an important step towards the rigorous understanding of the superheating field was made by Serfaty in \cite{Ser2}. There, a branch of vortexless local minimizers was constructed and shown to be unique, up to gauge transformations, in an explicit regime of applied fields extending far beyond the first critical field.

The three-dimensional homogeneous problem is considerably more involved. In two dimensions, vortices are localized near isolated points and can be described in terms of their topological degrees. In three dimensions, they concentrate instead along curves, and their analysis requires tools from geometric measure theory. Among the works that laid the foundations of the three-dimensional theory, we mention \cites{Riv,LinRiviere,sandier-surface-ball-construction,BalJerOrlSon1,BalJerOrlSon2}, where the limiting behavior of vortex filaments and the asymptotic structure of the Ginzburg--Landau energy were investigated in different regimes. More refined tools allowing one to describe vortex filaments at the $\ep$-level were later developed in \cite{roman-vortex-construction}; see also \cites{Rom-CMP,roman-sandier-serfaty,RSS2} for applications to the behavior of minimizers near the first critical field.

These $\ep$-level estimates were a key ingredient in extending the approach of Serfaty \cite{Ser2} to the three-dimensional homogeneous setting in \cite{Rom-CMP}. This extension established the persistence and uniqueness of the Meissner branch in suitable regimes of applied fields, providing the first three-dimensional counterpart of the two-dimensional theory. More recently, the approach was adapted to the two-dimensional Ginzburg--Landau model with pinning in \cite{diaz-roman-2d}.

The Ginzburg--Landau model with pinning has been extensively studied in two dimensions, in connection with a wide variety of questions concerning vortex nucleation, localization, and interaction. We refer to \cite{diaz-roman-2d} and the references therein for an overview of this literature. By contrast, its three-dimensional counterpart has received comparatively little mathematical attention. In \cite{DosSantos3D}, a non-magnetic problem with pinning and nontrivial Dirichlet boundary conditions was considered, and the limiting singular set was shown to form a geodesic link determined by the pinning landscape. Although most of their results concern the two-dimensional setting, Dos Santos, Rodiac, and Sandier \cite{DosRodSan} also considered the full Ginzburg--Landau energy in three dimensions with rapidly oscillating pinning terms. Assuming that the average of the pinning term in the periodic case, or its expectation in the random stationary case, is equal to $1$, they proved that, as $\ep\to0$, the limiting functional coincides with the one obtained in the homogeneous case by Baldo, Jerrard, Orlandi, and Soner \cite{BalJerOrlSon1}.

More recently, the development of a weighted counterpart of the $\ep$-level estimates introduced in \cite{roman-vortex-construction} was instrumental in the study of the first critical field for the three-dimensional magnetic Ginzburg--Landau functional with pinning in \cites{diaz-roman-3d,RomanPinningUpper25}. The tools developed in these works also provide a basis for addressing the different question considered here, concerning the persistence and uniqueness of the Meissner branch. They allow us to follow, in broad terms, the strategy used in the homogeneous setting, although its implementation in the presence of pinning is not direct. In particular, the construction and uniqueness arguments require elliptic estimates for equations whose coefficients depend on the pinning profile. Both the gain in integrability and the constants involved must be controlled uniformly in $\ep$, which requires a careful choice of gauges and suitable control of the oscillations of the weight.

\subsection{Main results}

Throughout the paper, $C>0$ denotes a constant independent of $\ep$ and $h_{\ex}$, whose value may change from line to line and may depend on the fixed parameters of the problem.

Since $\diver H_{0,\ex}=0$ in $\R^3$, in accordance with the nonexistence of magnetic monopoles in Maxwell's theory of electromagnetism, there exists a vector potential $A_{0,\ex}\in H^1_{\loc}(\R^3,\R^3)$ such that
$$
\curl A_{0,\ex}=H_{0,\ex}\quad\mathrm{in}\ \R^3.
$$
We define
$$
A_{\ex}\colonequals h_{\ex}A_{0,\ex}.
$$
Observe that we can choose $A_{0,\ex}$ in the Coulomb gauge in $\Omega$, namely,
$$
\left\{
\begin{array}{rcll}
	\diver A_{0,\ex}&=&0&\mathrm{in}\ \Omega\\
	A_{0,\ex}\cdot\nu&=&0&\mathrm{on}\ \partial\Omega,
\end{array}
\right.
$$
where hereafter $\nu$ denotes the outer unit normal to $\partial\Omega$.
The natural space for the minimization of $\glthreed$ in three dimensions is $H^1(\Omega,\C)\times[A_{\ex}+H_{\curl}]$, where
$$
H_{\curl}\colonequals\left\{A\in H^1_{\loc}(\R^3,\R^3)\ \middle|\  \curl A\in L^2(\R^3,\R^3)\right\}.
$$

Let us also introduce the homogeneous Sobolev space $\dot H^1(\R^3,\R^3)$, defined as the completion of $C_0^\infty(\R^3,\R^3)$ with respect to the norm $\|\nabla(\,\cdot\,)\|_{L^2(\R^3,\R^3)}$. We denote by
$$
\dot H^1_{\diver=0}\colonequals\left\{A\in\dot H^1(\R^3,\R^3) \ \middle|\  \diver A=0\ \mathrm{in}\ \R^3\right\}
$$
its divergence-free subspace. For every $A\in\dot H^1_{\diver=0}$, we have
\begin{equation}\label{space2}
	\|A\|_{\dot H^1_{\diver=0}}\colonequals\|\nabla A\|_{L^2(\R^3,\R^3)}=\|\curl A\|_{L^2(\R^3,\R^3)}.
\end{equation}

These spaces will be used below to formulate the convergence of the magnetic vector potential and to fix suitable gauges. Before stating our main results, we introduce a useful decomposition of the order parameter, whose underlying idea goes back to Lassoued and Mironescu \cite{lassoued-mironescu}, that isolates the effect of the pinning term in the energy. To this end, we first consider the scalar profile induced by the pinning term. More precisely, let $\pmin$ be the unique positive real-valued function satisfying
\begin{equation}\label{pde rho}
	\left\{
	\begin{array}{rcll}
		-\Delta \pmin&=&\dfrac{\pmin}{\ep^2}(a_\ep-\pmin^2)&\mathrm{in}\ \Omega\\
		\dfrac{\partial \pmin}{\partial \nu}&=&0&\mathrm{on}\ \partial\Omega.
	\end{array}
	\right.
\end{equation}
The function $\pmin$ represents the equilibrium modulus selected by the first two terms of \eqref{GLenergy} in the absence of a magnetic field. Indeed, while the potential term favors the pointwise relation $|\u|^2=a_\ep$, the gradient term penalizes the spatial variations required to reproduce $a_\ep$. Equation \eqref{pde rho} reflects the balance between these two competing effects, and $\pmin^2$ may therefore be interpreted as a regularized version of $a_\ep$.

Since $\pmin$ already incorporates the variations in the modulus of the order parameter induced by the pinning term, it is natural to factor out this profile by writing $\u=\pmin u$. The new order parameter $u$ then describes the remaining variations relative to the preferred scalar profile $\pmin$. Letting
$$
E_\ep(u,A)\colonequals\frac12\int_\Omega\left( |\nabla_Au|^2+\frac{1}{2\ep^2}(a_\ep-|u|^2)^2\right),\qquad E_\ep(u)=E_\ep(u,0),
$$
we have
\begin{equation}\label{lm mag decoupling}
	E_\ep(\pmin u,A)=E_\ep(\pmin)+\frac12\int_\Omega\left(\pmin^2|\nabla_Au|^2+\frac{\pmin^4}{2\ep^2}(1-|u|^2)^2\right).
\end{equation}
Therefore, the analysis of \eqref{GLenergy} can be reduced to that of the weighted Ginzburg--Landau functional
\begin{equation*}
	\frac12\int_\Omega\left(\pmin^2|\nabla_\A u|^2+\frac{\pmin^4}{2\ep^2}(1-|u|^2)^2\right)+\frac12\int_{\R^3}|\H-H_\ex|^2.
\end{equation*}
We refer to \cite[Section 2]{diaz-roman-3d} for details. The influence of the pinning term is thus incorporated into the spatially dependent weights $\pmin^2$ and $\pmin^4$. In particular, the estimates used in the construction and analysis of the Meissner solution, which we introduce next, must account for these weights and remain uniform as $\ep\to0$.

The Ginzburg--Landau model admits a distinguished configuration, unique modulo gauge transformations, that we refer to as the Meissner state, in reference to the Meissner effect in physics. It is obtained by minimizing the Ginzburg--Landau energy under the constraint $|\u|=\pmin$. More precisely, for each vector potential $A$, the phase $\phi_A$ is chosen to minimize
$$
\frac12\int_\Omega\pmin^2|A-\nabla\phi|^2.
$$
Equivalently, $A-\nabla\phi_A$ is in the weighted Coulomb gauge associated with the weight $\pmin^2$, namely,
\begin{equation}\label{eqphi0}
	\left\{
	\begin{array}{rcll}
		\diver(\pmin^2(A-\nabla\phi_A))&=&0&\mathrm{in}\ \Omega\\
		(A-\nabla\phi_A)\cdot\nu&=&0&\mathrm{on}\ \partial\Omega.
	\end{array}
	\right.
\end{equation}
This reduces the constrained problem to the minimization of a functional depending only on the vector potential; we refer to \cite[Section 2.3]{diaz-roman-3d} for the precise construction and discuss the weighted Coulomb gauge in more detail in Section~\ref{subsection:gauge-choices}.

We denote by $A_\ep^0$ the minimizer of this reduced functional, chosen so that
\begin{equation}\label{spaceforA0}
A_\ep^0-A_{0,\ex}\in\dot H^1_{\diver=0},
\end{equation}
and let $\phi_\ep^0\colonequals\phi_{A_\ep^0}$. The Meissner state then takes the form
$$
(\meisconf),
$$
where $\phi_\ep^0$ and $A_\ep^0$ depend on the domain $\Omega$, the normalized applied field $H_{0,\ex}$, and the weight function $\pmin$.

In view of \eqref{eqphi0}, the Helmholtz--Hodge decomposition yields the existence of a unique vector field $\Bfield$ such that
\begin{equation}\label{defBfield}
	\pmin^2(A_\ep^0-\nabla\phi_\ep^0)=\curl\Bfield\quad\mbox{in }\Omega,
\end{equation}
with
\begin{equation}\label{conditionsB0}
	\left\{
	\begin{array}{rcll}
		\diver\Bfield&=&0&\mbox{in }\Omega\\
		\Bfield\times\nu&=&0&\mbox{on }\partial\Omega.
	\end{array}
	\right.
\end{equation}
The estimates resulting directly from the minimization give, in particular,
$$
\|A_\ep^0-A_{0,\ex}\|_{\dot H^1_{\diver=0}}+\|\nabla\phi_\ep^0\|_{L^2(\Omega,\R^3)}+\|\curl\Bfield\|_{L^2(\Omega,\R^3)}\leq C.
$$

The Euler--Lagrange equation satisfied by $A_\ep^0$, which may be viewed as a weighted London equation for the vector potential associated with the induced magnetic field $\curl A_\ep^0$, is
$$
\curl^2(A_\ep^0-A_{0,\ex})+\chi_\Omega\curl\Bfield=0\quad\mbox{in }\R^3.
$$
From \eqref{spaceforA0}, we deduce that
$$
-\Delta(A_\ep^0-A_{0,\ex})=-\chi_\Omega\curl\Bfield\quad\mbox{in }\R^3.
$$
Since $\chi_\Omega\curl\Bfield\in L^2(\R^3,\R^3)$, interior elliptic regularity, applied with a cut-off supported in a sufficiently large ball containing $\overline\Omega$, yields
$$
\|A_\ep^0-A_{0,\ex}\|_{W^{2,2}(\Omega,\R^3)}\leq C.
$$
By Sobolev embedding, we then obtain
$$
\|A_\ep^0-A_{0,\ex}\|_{W^{1,6}(\Omega,\R^3)}\leq C.
$$
On the other hand, since $A_{0,\ex}$ is in the Coulomb gauge in $\Omega$, Lemma~\ref{lemma:div-curl} gives
$$
\|A_{0,\ex}\|_{W^{1,6}(\Omega,\R^3)}\leq C\|H_{0,\ex}\|_{L^6(\Omega,\R^3)}.
$$
Consequently,
$$
\|A_\ep^0\|_{W^{1,6}(\Omega,\R^3)}\leq C.
$$
By Sobolev embedding, we conclude that
\begin{equation}\label{estimateforA0}
	\|A_\ep^0\|_{L^\infty(\Omega,\R^3)}\leq C.
\end{equation}

Since $b\leq\pmin^2\leq1$ (see Proposition~\ref{prop:regularitypmin}), the classical result of Gr\"oger \cite{Groger} (see also \cite{Meyers} for Dirichlet boundary conditions) establishes that there exists $\pzero>2$ such that
\begin{equation}\label{estimateforgradphip0}
	\|\nabla\phi_\ep^0\|_{L^{\pzero}(\Omega,\R^3)}\leq C\|A_\ep^0\|_{L^{\pzero}(\Omega,\R^3)},
\end{equation}
where both $\pzero$ and $C$ depend on $b$ but not on $\ep$. For our analysis, however, we need \eqref{estimateforgradphip0} to hold uniformly in $\ep$ for the exponents specified in the statements of our main results. This is one of the points where the presence of the spatially dependent weight creates substantial additional difficulties. In particular, both the gain in integrability and the constants in the relevant elliptic estimates must be controlled uniformly in $\ep$.

In the course of writing this article, we realized that there is a mistake in the proof of \cite[Proposition 2.2]{diaz-roman-3d} concerning the regularity of $\nabla\phi_\ep^0$. There, it was incorrectly concluded that \eqref{estimateforgradphip0} holds uniformly in $\ep$ for any $\pzero>2$. In fact, extending this estimate beyond the exponent provided by the general theory requires either further assumptions on $\pmin$ or an $\ep$-dependent constant which, loosely speaking, depends on a suitable control of the oscillations of $\pmin^2$ at sufficiently small scales; see, for instance, \cite{SunWang}.In a work in progress \cite{DRY}, we will show that, under the main assumption on $\pmin$ in \cite{diaz-roman-3d}, estimate \eqref{estimateforgradphip0} holds for any $\pzero>2$, with a constant $C=C(\ep)=\lep^M$ for some $M>0$. This estimate still suffices to establish the results of that paper.

To ensure that \eqref{estimateforgradphip0} holds for the required exponent, we fix $\pzero>2$, as specified in the statements of our main results, and assume hereafter that there exist constants $r_0>0$ and $\de_0=\de_0(b,\pzero,r_0)>0$ such that, for all sufficiently small $\ep$,
\begin{equation}\label{conditionBMO}
	[\pmin^2]_{BMO_{r_0}(\Omega)}\leq\frac{\de_0}{2}.
\end{equation}
Under this assumption, \cite[Theorem 2]{SunWang} yields \eqref{estimateforgradphip0}. Let us also remark that \eqref{conditionBMO} is automatically satisfied whenever $b$ is sufficiently close to $1$.

In view of \eqref{defBfield}, combining \eqref{estimateforA0} and \eqref{estimateforgradphip0}, we also obtain
\begin{equation}\label{estimateforB0}
	\|\curl\Bfield\|_{L^{\pzero}(\Omega,\R^3)}\leq C.
\end{equation}
Moreover, in view of \eqref{conditionsB0}, \cite[Corollary 3.4]{amrouche-seloula} yields
\begin{equation}\label{estimateforB0W1p}
	\|\Bfield\|_{W^{1,\pzero}(\Omega,\R^3)}\leq C.
\end{equation}

We are now in a position to state our first main result. It establishes the existence of a locally minimizing vortexless configuration in a range of applied fields extending far beyond the first critical field. In addition to its local minimality, we obtain quantitative estimates showing that this configuration is close, up to gauge transformations, to the Meissner state introduced above.

\begin{theorem}\label{theorem:meissner}
	Let $\alpha\in\left(0,\frac14\right)$ and suppose that
	$$
	\pzero>\frac3{1-3\alpha}.
	$$
	There exists $\ep_0>0$ such that, for any $\ep<\ep_0$, if $h_{\ex}\leq\ep^{-\alpha}$, then there exists a vortexless configuration
	$(\u,\A)=(\meisdecom)$ which locally minimizes $\glthreed$ in $H^1(\Omega,\C)\times[A_{\ex}+H_{\curl}]$ and satisfies the following properties as $\ep\to0$:
	\begin{enumerate}[label=\textnormal{\bfseries (\alph*)},ref=(\alph*),leftmargin=0pt,labelsep=*,itemindent=*,itemsep=5pt,topsep=5pt]
		\item\label{part-a} $\glthreed(\meisconf)+o(1)\leq\glthreed(\u,\A)\leq\glthreed(\meisconf)$.
		\item\label{part-b} $\|1-|u|\|_{L^\infty(\Omega)}=o(1)$ and $\||\u|-\pmin\|_{L^\infty(\Omega)}=o(1)$.
		\item\label{part-c} There exists $\de>0$ such that
		$$
		\fenthreedOmega(u,A)\leq C\ep^{1+\de}.
		$$
		\item\label{part-d} $\|\mu(u,A)\|_{C_T^{0,1}(\Omega,\R^3)^*}=o(1)$ and $\|\mu(\pmin^{-1}\u,\A)\|_{C_T^{0,1}(\Omega,\R^3)^*}=o(1)$.
		\item\label{part-e} Up to a gauge transformation, we have
		$$
		\inf_{\theta\in[0,2\pi]}\|u-e^{i\theta}\|_{H^1(\Omega,\C)}+\|A\|_{\dot H^1_{\diver=0}}=o(1).
		$$
		\item\label{part-f} Up to a gauge transformation, $(\pmin^{-1}\u,\A)$ converges to $(\phase,h_\ex A_\ep^0)$. More precisely, we have
		$$
		\inf_{\theta\in[0,2\pi]}\|\pmin^{-1}\u-e^{i\theta}\phase\|_{H^1(\Omega,\C)}+\|\A-h_{\ex}A_\ep^0\|_{\dot H^1_{\diver=0}}=o(1).
		$$
	\end{enumerate}
\end{theorem}

\begin{remark}
	The condition $\pzero>\frac3{1-3\alpha}$ approaches the natural threshold $\pzero>3$ as $\alpha\to0$. The requirement $\pzero>3$ is related to the H\"older regularity needed to estimate the vorticity term. More precisely, the argument uses
	$$
	W^{1,\pzero}(\Omega)\hookrightarrow C^{0,\beta_0}(\overline\Omega),
	\qquad
	\beta_0=1-\frac3{\pzero},
	$$
	together with a quantitative vorticity estimate that provides a gain of order $\ep^{\frac{\beta_0}{3}}$. The condition $\frac{\beta_0}{3}>\alpha$ is exactly equivalent to the condition on $\pzero$ above.
\end{remark}

Our last result concerns the uniqueness, up to a gauge transformation, of locally minimizing configurations satisfying a suitable free-energy bound after factoring out the Meissner state. The vortexless character of the configuration is not assumed, but instead follows from the hypotheses. In particular, this result applies to the locally minimizing branch constructed in Theorem~\ref{theorem:meissner}.

\begin{theorem}\label{theorem:uniqueness}
	Let $\alpha\in(0,1)$ and suppose that
	$$
	\pzero>\frac3{1-\alpha}.
	$$
	There exist constants $c_0>0$ and $\ep_0>0$ such that, for any $\ep<\ep_0$, if $h_{\ex}\leq\ep^{-\alpha}$, then a configuration $(\u,\A)=(\meisdecom)$ which locally minimizes $\glthreed$ in $H^1(\Omega,\C)\times[A_{\ex}+H_{\curl}]$ and satisfies
	$\fenthreedOmega(u,A)\leq c_0\ep$ is unique up to a gauge transformation.
\end{theorem}

\begin{remark}
	Theorem~\ref{theorem:uniqueness} applies, in particular, to the locally minimizing configuration constructed in Theorem~\ref{theorem:meissner}. Indeed, for some $\de>0$,
	$$
	\fenthreedOmega(u,A)\leq C\ep^{1+\de}=o(\ep).
	$$
	Consequently, for any fixed $c_0>0$, we have
	$$
	\fenthreedOmega(u,A)\leq c_0\ep
	$$
	for all sufficiently small $\ep$. Moreover, in the range $\alpha\in(0,\frac14)$ covered by Theorem~\ref{theorem:meissner},
	$$
	\frac3{1-3\alpha}>\frac3{1-\alpha},
	$$
	so the condition on $\pzero$ required for existence automatically implies the one required for uniqueness. Therefore, the locally minimizing branch constructed in Theorem~\ref{theorem:meissner} is unique up to a gauge transformation without any additional assumption on $\pzero$.
	
	For this common range of existence and uniqueness, the required exponent approaches the threshold $\pzero>3$ as $\alpha\to0$, which is natural in view of the H\"older embedding
	$$
	W^{1,\pzero}(\Omega)\hookrightarrow C^{0,1-\frac3{\pzero}}(\overline\Omega).
	$$
	On the other hand, the required exponent approaches $12$ as $\alpha\to\frac14$.
\end{remark}

Together, Theorems~\ref{theorem:meissner} and~\ref{theorem:uniqueness} provide a rigorous description of the persistence of the Meissner branch in the three-dimensional pinned setting. They show that locally minimizing vortexless configurations continue to exist, and remain unique within a natural quantitatively controlled class, for applied fields extending far beyond the first critical field. In particular, the loss of global minimality at $H_{c_1}$ does not correspond to the disappearance of the Meissner branch. This is consistent with the expected superheating scenario, according to which this branch should remain locally stable until a much larger applied field, of order $\ep^{-1}$.

These results constitute, to the best of our knowledge, the first rigorous analysis of the persistence and uniqueness of the Meissner branch for the three-dimensional magnetic Ginzburg--Landau functional with a spatially dependent pinning term. Although the ranges obtained here do not reach the expected scale of the superheating field, they provide an algebraic regime of locally stable vortexless solutions and extend the variational picture previously established for the homogeneous model to the pinned setting.

\subsection*{Acknowledgments}

This work was partially funded by ANID FONDECYT 1231593.

\subsection*{AI declaration}

The mathematical content of this paper, including its results and proof strategies, was developed by the author. Microsoft 365 Copilot, based on the GPT-5 reasoning model, was used to assist with the revision of the manuscript, including language editing, improvements in clarity and organization, and checks of notation, consistency, and selected mathematical calculations. The author carefully reviewed and independently verified all AI-assisted changes and takes full responsibility for the content of the manuscript.

%%%%%%%%%%%%%%%%%%%%%%%%%%%%%%%%%%%%%%%%%%%%%%%%%
\section{Preliminaries}
%%%%%%%%%%%%%%%%%%%%%%%%%%%%%%%%%%%%%%%%%%%%%%%%%

In this section, we collect some preliminary results that will be used throughout the paper. We first introduce the different gauge choices needed in our analysis and recall some classical estimates. We then establish some basic properties of the pinning profile $\pmin$ and of critical points of the Ginzburg--Landau energy. Finally, we recall a useful splitting of the energy and prove the clearing-out and vorticity estimates that will play an important role in the proofs of our main results.

\subsection{Gauge choices}\label{subsection:gauge-choices}

We begin by introducing the three gauge choices that will be used throughout the paper. The first two are defined in bounded domains, while the third one is a global divergence-free gauge in $\R^3$.

\begin{definition}
	We say that a vector field $A$ is in the Coulomb gauge in a bounded domain $\Theta$ if
	\begin{equation}\label{Coulombgauge}
		\left\{ 
		\begin{array}{rcll}
			\diver A&=&0&\mathrm{in}\ \Theta\\
			A\cdot \nu&=&0&\mathrm{on}\ \partial\Theta.
		\end{array}\right.
	\end{equation}
\end{definition}

Every vector field $A\in L^2(\Theta,\R^3)$ can be put in this gauge. Indeed, there exists a function $\phi\in H^1(\Theta)$, unique up to an additive constant, characterized by
$$
\int_\Theta\nabla\phi\cdot\nabla\psi=\int_\Theta A\cdot\nabla\psi,\quad\forall\psi\in H^1(\Theta).
$$
Equivalently, $\phi$ is the weak solution of
\begin{equation*}
	\left\{
	\begin{array}{rcll}
		\Delta\phi&=&\diver A&\mathrm{in}\ \Theta\\
		\dfrac{\partial\phi}{\partial\nu}&=&A\cdot\nu&\mathrm{on}\ \partial\Theta.
	\end{array}
	\right.
\end{equation*}
It follows that $A-\nabla\phi$ satisfies \eqref{Coulombgauge}. After extending $\phi$ to a function in $H^1_{\loc}(\R^3)$, the corresponding gauge transformation is defined in the whole space.

\smallskip
The presence of the pinning profile $\pmin$ leads naturally to the following weighted version of the Coulomb gauge, which arises in the construction of the Meissner state and in the estimates for the phase.

\begin{definition}
	We say that $A$ is in the weighted Coulomb gauge associated with the weight $\pmin^2$ if
	\begin{equation}\label{weightedCoulombgauge}
		\left\{
		\begin{array}{rcll}
			\diver(\pmin^2A)&=&0&\mathrm{in}\ \Omega\\
			A\cdot\nu&=&0&\mathrm{on}\ \partial\Omega.
		\end{array}
		\right.
	\end{equation}
\end{definition}

As in the previous case, every vector field $A\in L^2(\Omega,\R^3)$ can be put in this gauge. Indeed, since $b\leq\pmin^2\leq1$, there exists a function $\phi\in H^1(\Omega)$, unique up to an additive constant, characterized by
$$
\int_\Omega\pmin^2\nabla\phi\cdot\nabla\psi=\int_\Omega\pmin^2A\cdot\nabla\psi,\quad\forall\psi\in H^1(\Omega).
$$
Equivalently, $\phi$ is the weak solution of
\begin{equation}\label{weightedCoulombgaugephi}
	\left\{
	\begin{array}{rcll}
		\diver(\pmin^2\nabla\phi)&=&\diver(\pmin^2A)&\mathrm{in}\ \Omega\\
		\dfrac{\partial\phi}{\partial\nu}&=&A\cdot\nu&\mathrm{on}\ \partial\Omega.
	\end{array}
	\right.
\end{equation}
It follows that $A-\nabla\phi$ satisfies \eqref{weightedCoulombgauge}. As before, $\phi$ may be extended to a function in $H^1_{\loc}(\R^3)$, so that the corresponding gauge transformation is defined in the whole space.

\begin{remark}
	As discussed in the introduction, the regularity of the solution $\phi$ of \eqref{weightedCoulombgaugephi} depends on the regularity and oscillation of the weight $\pmin^2$. Under assumption \eqref{conditionBMO}, if $A\in L^{\pzero}(\Omega,\R^3)$, then \cite[Theorem 2]{SunWang} yields
	$$
	\|\nabla\phi\|_{L^{\pzero}(\Omega,\R^3)}\leq C\|A\|_{L^{\pzero}(\Omega,\R^3)},
	$$
	where $C$ is independent of $\ep$.
\end{remark}

For the global gauge choice, we recall some additional properties of the homogeneous Sobolev space introduced in the previous section. By Sobolev embedding, there exists a constant $C>0$ such that
\begin{equation}\label{space1}
	\|A\|_{L^6(\R^3,\R^3)}\leq C\|\nabla A\|_{L^2(\R^3,\R^3)}
\end{equation}
for every $A\in\dot H^1(\R^3,\R^3)$. Moreover, by \cite[Proposition 2.4]{KozSoh},
$$
\dot H^1(\R^3,\R^3)=\left\{A\in L^6(\R^3,\R^3)\ \middle|\ \nabla A\in L^2(\R^3,\R^3)\right\}.
$$
In particular, combining \eqref{space1} and \eqref{space2}, for every $A\in\dot H^1_{\diver=0}$ we have
$$
\|A\|_{L^6(\R^3,\R^3)}\leq C\|\curl A\|_{L^2(\R^3,\R^3)}.
$$

We finally introduce the global divergence-free gauge.

\begin{definition}
	We say that a vector field $A$ is in the global divergence-free gauge if
	$$
	A\in\dot H^1_{\diver=0}.
	$$
\end{definition}

The fact that every configuration can be put in this gauge follows from the decomposition established in \cite[Lemma 2.1]{Rom-CMP}.

\begin{lemmaA}\label{lemma:decompositionR3}
	Every vector field $A\in H_{\curl}$ admits a decomposition
	$$
	A=\curl B+\nabla\Phi,
	$$
	where $B,\curl B\in\dot H^1_{\diver=0}$ and $\Phi\in H^2_{\loc}(\R^3)$.
\end{lemmaA}

Indeed, after the gauge transformation 
$$
u\to ue^{-i\Phi},\qquad A\to A-\nabla\Phi,
$$
we obtain $A=\curl B\in\dot H^1_{\diver=0}$.

\subsection{Some classical estimates}

We will repeatedly use the following classical estimates for vector fields in the Coulomb gauge, which follow directly from \cite[Corollaries 3.4 and 3.5]{amrouche-seloula}.

\begin{lemmaA}\label{lemma:div-curl}
	Let $\Theta$ be a bounded and simply connected domain of $\R^3$ with $C^2$ boundary, and let $A$ satisfy \eqref{Coulombgauge} in $\Theta$. Then, for any $p\in(1,\infty)$,
	\begin{align*}
		\|A\|_{W^{1,p}(\Theta,\R^3)}&\leq C\|\curl A\|_{L^p(\Theta,\R^3)}\\
		\|A\|_{W^{2,p}(\Theta,\R^3)}&\leq C\|\curl A\|_{W^{1,p}(\Theta,\R^3)}.
	\end{align*}
	Here, $C$ depends only on $\Theta$ and $p$.
\end{lemmaA}

We will also use the following estimate, which is the Neumann counterpart of \cite[Lemma A.2]{BBH-CVPDE}. The proof follows the same argument, with only minor modifications.

\begin{lemmaA}\label{lemma:BBH}
	Let $\Theta$ be a bounded domain of $\R^N$, with $N\geq2$, having $C^2$ boundary. If $v$ satisfies
	$$
	\left\lbrace 
	\begin{array}{rcll}
		-\Delta v&=&f&\mathrm{in}\ \Theta\\
		\frac{\partial v}{\partial \nu}&=&0&\mathrm{on}\ \partial\Theta,
	\end{array}
	\right.
	$$
	where $f$ is such that $\int_\Theta f=0$, then there exists a constant $C>0$ depending only on $\Theta$ such that
	$$
	\|\nabla v\|_{L^\infty(\Theta)}^2\leq C\|v\|_{L^\infty(\Theta)}\left(\|\Delta v\|_{L^\infty(\Theta)}+\|v\|_{L^\infty(\Theta)}\right).
	$$
\end{lemmaA}

\subsection{Properties of the pinning profile}

We next establish some basic properties of the pinning profile $\pmin$ defined by \eqref{pde rho}.

\begin{prop}\label{prop:regularitypmin}
	We have $b\leq\pmin^2\leq1$ and
	$$
	\|\nabla\pmin\|_{L^\infty(\Omega,\R^3)}\leq C\ep^{-1}.
	$$
\end{prop}

\begin{proof}
	Testing \eqref{pde rho} against $(\pmin-1)_+\in H^1(\Omega)$, we find
	$$
	\int_{\{\pmin>1\}}|\nabla\pmin|^2=\frac1{\ep^2}\int_{\{\pmin>1\}}(\pmin-1)\pmin(a_\ep-\pmin^2)\leq0,
	$$
	where the inequality follows from $a_\ep-\pmin^2<0$ when $\pmin>1$. From this we deduce that $|\{\pmin>1\}|=0$, and hence $\pmin\leq1$.
	
	Analogously, by testing \eqref{pde rho} against $(\pmin-\sqrt{b})_-=-\min\left\{\pmin-\sqrt{b},0\right\}$, we deduce that $\sqrt{b}\leq\pmin$. Therefore,
	$$
	b\leq\pmin^2\leq1.
	$$
	
	On the other hand, Lemma~\ref{lemma:BBH}, together with the bounds $b\leq\pmin^2\leq1$ and $b\leq a_\ep\leq1$, yields
	\begin{align*}
		\|\nabla\pmin\|_{L^\infty(\Omega,\R^3)}^2&\leq C\|\pmin\|_{L^\infty(\Omega)}\left(\|\Delta\pmin\|_{L^\infty(\Omega)}+\|\pmin\|_{L^\infty(\Omega)}\right)\\
		&\leq C\left\|\dfrac{\pmin}{\ep^2}(a_\ep-\pmin^2)\right\|_{L^\infty(\Omega)}+C\\
		&\leq\frac{C}{\ep^2}.
	\end{align*}
	The result follows.
\end{proof}

\subsection{Critical points of the Ginzburg--Landau energy}

We now collect some properties of critical points of $\glthreed$ that will be used throughout the paper.

\begin{prop}[Ginzburg--Landau equations]
	If $(\u,\A)\in H^1(\Omega,\C)\times[A_{\ex}+H_{\curl}]$ is a critical point of $\glthreed$, then $(\u,\A)$ satisfies the system of equations
	\begin{equation}\label{GLeq}
		\left\lbrace 
		\begin{array}{rcll}
			-(\nabla_\A)^2\u&=&\displaystyle\frac1{\ep^2}\u(a_\ep-|\u|^2)&\mathrm{in}\ \Omega\\
			\curl(\H-H_\ex)&=&(i\u,\nabla_\A\u)\chi_\Omega&\mathrm{in}\ \R^3\\
			\nabla_\A\u\cdot\nu&=&0&\mathrm{on}\ \partial\Omega\\
			\lbrack\H-H_\ex\rbrack\times\nu&=&0&\mathrm{on}\ \partial\Omega.
		\end{array}
		\right.
	\end{equation}
	Here, $\chi_\Omega$ is the characteristic function of $\Omega$, $[\,\cdot\,]$ denotes the jump across $\partial\Omega$,
	$$
	\nabla_\A\u\cdot\nu=\sum_{j=1}^3(\partial_j\u-i\A_j\u)\nu_j,
	$$
	and the covariant Laplacian $(\nabla_\A)^2$ is defined by
	$$
	(\nabla_\A)^2\u=(\diver-i\A\cdot)\nabla_\A\u.
	$$
\end{prop}
\begin{proof}
	The proof follows by standard variations of $\glthreed$ with respect to $\u$ and $\A$ and is analogous to the proof in the homogeneous case given in \cite[Proposition 2.1]{Rom-CMP}.
\end{proof}

\begin{prop}\label{prop:regularity}
	Let $(\u,\A)$ be a critical point of $\glthreed$. Then
	\begin{equation}\label{boundu}
		\|\u\|_{L^\infty(\Omega,\C)}\leq1\quad\mbox{and}\quad\left\|\frac{\u}{\pmin}\right\|_{L^\infty(\Omega,\C)}\leq1.
	\end{equation}
	Moreover, if $\A$ is in the Coulomb gauge in $\Omega$, then
	\begin{equation}\label{boundnablau}
		\|\nabla\u\|_{L^\infty(\Omega,\C^3)}\leq C\left(\frac1\ep+\|\A\|_{L^\infty(\Omega,\R^3)}\right).
	\end{equation}
\end{prop}

\begin{proof}
	By taking the inner product of the first equation in \eqref{GLeq} with $\u$ and using the first boundary condition, a direct computation shows that $|\u|$ satisfies
	$$
	\left\lbrace 
	\begin{array}{rcll}
		-\frac12\Delta|\u|^2&=&\displaystyle\frac1{\ep^2}|\u|^2(a_\ep-|\u|^2)-|\nabla_\A\u|^2&\mathrm{in}\ \Omega\\
		\nabla|\u|^2\cdot\nu&=&0&\mathrm{on}\ \partial\Omega.
	\end{array}
	\right.
	$$
	Testing against $(|\u|^2-1)_+\in H^1(\Omega)$, we find
	$$
	\frac12\int_{\{|\u|^2>1\}}|\nabla|\u|^2|^2=\frac1{\ep^2}\int_{\{|\u|^2>1\}}|\u|^2(a_\ep-|\u|^2)(|\u|^2-1)-|\nabla_\A\u|^2(|\u|^2-1)\leq0,
	$$
	where the inequality follows from $a_\ep-|\u|^2<0$ when $|\u|^2>1$. From this we deduce that $|\{|\u|^2>1\}|=0$, and hence $|\u|^2\leq1$ in $\Omega$.
	
	Now, writing $u=\pmin^{-1}\u$, by combining the first equation and boundary condition in \eqref{GLeq} with \eqref{pde rho}, we find
	$$
	\left\lbrace 
	\begin{array}{rcll}
		-(\diver-i\A\cdot)(\pmin^2\nabla_\A u)&=&\displaystyle\frac1{\ep^2}\pmin^4u(1-|u|^2)&\mathrm{in}\ \Omega\\
		\nabla_\A u\cdot\nu&=&0&\mathrm{on}\ \partial\Omega.
	\end{array}
	\right.
	$$
	Taking the inner product with $u$, a direct computation shows that $|u|$ satisfies
	$$
	\left\lbrace 
	\begin{array}{rcll}
		-\frac12\diver(\pmin^2\nabla|u|^2)&=&\displaystyle\frac{\pmin^4}{\ep^2}|u|^2(1-|u|^2)-\pmin^2|\nabla_\A u|^2&\mathrm{in}\ \Omega\\
		\nabla|u|^2\cdot\nu&=&0&\mathrm{on}\ \partial\Omega.
	\end{array}
	\right.
	$$
	Testing against $(|u|^2-1)_+\in H^1(\Omega)$, we find
	$$
	\frac12\int_{\{|u|^2>1\}}\pmin^2|\nabla|u|^2|^2=\frac1{\ep^2}\int_{\{|u|^2>1\}}\pmin^4|u|^2(1-|u|^2)(|u|^2-1)-\pmin^2|\nabla_\A u|^2(|u|^2-1)\leq0,
	$$
	where the inequality follows from $1-|u|^2<0$ when $|u|^2>1$. From this we deduce that $|\{|u|^2>1\}|=0$, and hence $|u|^2\leq1$ in $\Omega$.
	
	On the other hand, if $\A$ is in the Coulomb gauge in $\Omega$, by expanding the covariant Laplacian, we find
	$$
	-\Delta\u=\frac1{\ep^2}\u(a_\ep-|\u|^2)-2i\A\cdot\nabla\u-|\A|^2\u,
	$$
	together with the boundary condition $\nabla\u\cdot\nu=0$ on $\partial\Omega$. Hence, by applying Lemma~\ref{lemma:BBH} and using $\|\u\|_{L^\infty(\Omega,\C)}\leq1$, we find
	$$
	\|\nabla\u\|_{L^\infty(\Omega,\C^3)}^2\leq C\left(\frac1{\ep^2}+\|\A\|_{L^\infty(\Omega,\R^3)}\|\nabla\u\|_{L^\infty(\Omega,\C^3)}+\|\A\|_{L^\infty(\Omega,\R^3)}^2\right),
	$$
	from where it follows that
	$$
	\|\nabla\u\|_{L^\infty(\Omega,\C^3)}\leq C\left(\frac1\ep+\|\A\|_{L^\infty(\Omega,\R^3)}\right).
	$$
\end{proof}

\begin{prop}\label{prop:regularityA}
	Let $(\u,\A)$ be a critical point of $\glthreed$. If $\A$ is in the Coulomb gauge in $\Omega$, then
	$$
	\|\A\|_{L^\infty(\Omega,\R^3)}\leq C\left(\left\|\nabla_\A\left(\frac{\u}{\pmin}\right)\right\|_{L^2(\Omega,\C^3)}+h_\ex\right).
	$$
\end{prop}

\begin{proof}
	Since $(\u,\A)$ is a critical point of $\glthreed$, it satisfies the second Ginzburg--Landau equation
	\begin{equation}\label{secondGLeq}
		\curl^2(\A-A_\ex)=(i\u,\nabla_\A\u)\chi_\Omega\quad\mathrm{in}\ \R^3.
	\end{equation}
	Since $\A-A_\ex\in H_{\curl}$, we have
	$$
	\curl(\A-A_\ex)\in L^2(\R^3,\R^3).
	$$
	Moreover,
	$$
	\diver\curl(\A-A_\ex)=0\quad\mathrm{in}\ \R^3,
	$$
	and, by \eqref{secondGLeq}, \eqref{boundu}, and the fact that $\nabla_\A\u\in L^2(\Omega,\C^3)$,
	$$
	\curl^2(\A-A_\ex)\in L^2(\R^3,\R^3).
	$$
	Using the classical whole-space div--curl identity
	$$
	\|\nabla\curl(\A-A_\ex)\|_{L^2(\R^3,\R^3)}^2
	=
	\|\curl^2(\A-A_\ex)\|_{L^2(\R^3,\R^3)}^2
	+
	\|\diver\curl(\A-A_\ex)\|_{L^2(\R^3)}^2,
	$$
	we obtain
	$$
	\curl(\A-A_\ex)\in H^1(\R^3,\R^3)
	$$
	and
	$$
	\|\nabla\curl(\A-A_\ex)\|_{L^2(\R^3,\R^3)}
	=
	\|\curl^2(\A-A_\ex)\|_{L^2(\R^3,\R^3)}.
	$$
	Since
	$$
	\curl(\A-A_\ex)\in H^1(\R^3,\R^3)\subset\dot H^1(\R^3,\R^3)
	$$
	and this vector field is divergence free, we also have
	$$
	\curl(\A-A_\ex)\in\dot H^1_{\diver=0}.
	$$
	
	Writing $u=\pmin^{-1}\u$, we observe that
	\begin{multline*}
		(i\u,\nabla_\A\u)=(i\pmin u,\nabla_\A(\pmin u))=\pmin(iu,u\nabla\pmin+\pmin\nabla_\A u)\\
		=\pmin\nabla\pmin(iu,u)+\pmin^2(iu,\nabla_\A u)=\pmin^2(iu,\nabla_\A u).
	\end{multline*}
	Using \eqref{secondGLeq}, \eqref{boundu}, and the fact that $\pmin\leq1$, we obtain
	$$
	\|\nabla\curl(\A-A_\ex)\|_{L^2(\R^3,\R^3)}=\|\curl^2(\A-A_\ex)\|_{L^2(\R^3,\R^3)}\leq\|\nabla_\A u\|_{L^2(\Omega,\C^3)}.
	$$
	In view of \eqref{space1}, we then find that
	$$
	\|\curl(\A-A_\ex)\|_{L^6(\R^3,\R^3)}\leq C\|\nabla_\A u\|_{L^2(\Omega,\C^3)}.
	$$
	Using Lemma~\ref{lemma:div-curl}, and recalling that both $\A$ and $A_\ex$ are in the Coulomb gauge in $\Omega$, we conclude that
	$$
	\|\A-A_\ex\|_{W^{1,6}(\Omega,\R^3)}
	\leq C\|\curl(\A-A_\ex)\|_{L^6(\Omega,\R^3)}
	$$
	and
	$$
	\|A_\ex\|_{W^{1,6}(\Omega,\R^3)}
	\leq C\|\curl A_\ex\|_{L^6(\Omega,\R^3)}.
	$$
	We then deduce that
	$$
	\|\A\|_{W^{1,6}(\Omega,\R^3)}
	\leq C\left(\|\nabla_\A u\|_{L^2(\Omega,\C^3)}+\|\curl A_\ex\|_{L^6(\Omega,\R^3)}\right).
	$$
	Finally, by Sobolev embedding and recalling that $\curl A_\ex=h_\ex H_{0,\ex}$, we find
	$$
	\|\A\|_{L^\infty(\Omega,\R^3)}
	\leq C\left(\|\nabla_\A u\|_{L^2(\Omega,\C^3)}+h_\ex\right).
	$$
	This concludes the proof.
\end{proof}

\subsection{Energy and vorticity estimates}

We conclude this section with three ingredients that will play a central role in the proofs of our main results. We first recall the energy splitting established in \cite[Proposition 1.1]{diaz-roman-3d}, and then prove a clearing-out estimate and a quantitative estimate for the vorticity.

\begin{propA}[Energy splitting]
	Given any configuration $(\u,\A)\in\minsp$, letting $(u,A)$ be defined via the relation
	$$
	(\u,\A)=(\pmin u\phase,A+h_\ex A_\ep^0),
	$$
	we have
	\begin{equation}\label{Split3D}
		\glthreed(\u,\A)=\glthreed(\pmin\phase,h_\ex A_\ep^0)+\fenweightthreed(u,A)-h_\ex\int_\Omega\mu(u,A)\cdot\Bfield+R_0,
	\end{equation}
	where
	$$
	\fenweightthreed(u,A)=\frac12\int_\Omega\left(\pmin^2|\nabla_Au|^2+\frac{\pmin^4}{2\ep^2}(1-|u|^2)^2\right)+\frac12\int_{\R^3}|\curl A|^2
	$$
	and
	$$
	R_0\colonequals\frac{h_\ex^2}{2}\int_\Omega\frac{|\curl\Bfield|^2}{\pmin^2}(|u|^2-1).
	$$
\end{propA}

\begin{lemma}[Clearing-out estimate]\label{lemma:clearingout}
	Let $u\in H^1(\Omega,\C)$ satisfy $\|\nabla|u|\|_{L^\infty(\Omega,\R^3)}\leq C\ep^{-1}$. Then
	$$
	(1-|u|(x_0))^5\leq\frac{C}{\ep^3}\int_{B_\ep(x_0)\cap\Omega}(1-|u|^2)^2,\quad\forall x_0\in\overline\Omega.
	$$
\end{lemma}

\begin{proof}
	The proof is an adaptation of \cite[Lemma III.3]{BetBreOrl}. Let
	$$
	k\colonequals1-|u|(x_0).
	$$
	We may assume that $k>0$, since otherwise the result trivially holds. For any $x\in\Omega$, we have
	$$
	\big||u|(x)-|u|(x_0)\big|\leq\|\nabla|u|\|_{L^\infty(\Omega,\R^3)}|x-x_0|\leq\frac{C}{\ep}|x-x_0|.
	$$
	Letting $r\colonequals\frac{k\ep}{2C}$, we obtain
	$$
	1-|u|(x)\geq k-\frac{C}{\ep}|x-x_0|\geq\frac{k}{2}\quad\forall x\in B_r(x_0)\cap\Omega.
	$$
	Since $0<k\leq1$, by increasing $C$ if necessary, we have $r<\ep$, and hence $B_r(x_0)\subset B_\ep(x_0)$. 
	
	The $C^2$ regularity of $\partial\Omega$ implies that there exists a constant $c>0$, depending only on $\Omega$, such that
	$$
	|B_r(x_0)\cap\Omega|\geq cr^3.
	$$
	Hence
	\begin{multline*}
		\int_{B_\ep(x_0)\cap\Omega}(1-|u|^2)^2\geq\int_{B_r(x_0)\cap\Omega}(1-|u|^2)^2=\int_{B_r(x_0)\cap\Omega}(1-|u|)^2(1+|u|)^2\\
		\geq\int_{B_r(x_0)\cap\Omega}(1-|u|)^2\geq c\frac{k^2}{4}r^3=\frac{c}{32C^3}k^5\ep^3.
	\end{multline*}
	The inequality is thus proved.
\end{proof}

\begin{prop}[Quantitative vorticity estimate]\label{prop:vorticityestimate}
	Let $(u,A)\in H^1(\Omega,\C)\times H^1(\Omega,\R^3)$ be such that $\fenthreedOmega(u,A)\leq\ep^\lambda$, for some $\lambda\in(0,1]$, and let $\beta\in(0,1]$.
	
	If $\lambda>\frac13$, there exists $C_\beta>0$ depending on $\beta$ and $\Omega$ such that
	$$
	\|\mu(u,A)\|_{C_T^{0,\beta}(\Omega,\R^3)^*}\leq C_\beta\ep^{\frac{\beta}{3}}\fenthreedOmega(u,A).
	$$
	
	If $\lambda\leq\frac13$, then, for any $\gamma\in(0,\lambda)$, there exists $C_\beta>0$ depending on $\beta$ and $\Omega$ such that
	$$
	\|\mu(u,A)\|_{C_T^{0,\beta}(\Omega,\R^3)^*}\leq C_\beta\ep^{\beta(\lambda-\gamma)}\fenthreedOmega(u,A).
	$$
	
	The estimates hold for all sufficiently small $\ep$, where the corresponding threshold may depend on $\lambda$ and, in the second case, on $\gamma$.
\end{prop}

\begin{proof}
	We consider, for $\ep$ sufficiently small and $\de=\de(\ep)=\ep^\sigma$, with $\sigma\leq\frac{1+\lambda}{4}$ to be fixed later, the grid $\GG(b_\ep,\de)$ associated with $(u,A)$ by \cite[Lemma 2.1]{roman-vortex-construction}. In particular, using the same notation as in that lemma, we have
	\begin{align}
		|u|&>5/8\quad\mathrm{on}\ \RR_1(\GG(b_\ep,\de)),\notag\\
		I_\ep^1\colonequals\int\limits_{\RR_1(\GG(b_\ep,\de))}e_\ep(u,A)d\HH^1&\leq C\de^{-2}\fenthreedOmega(u,A)\leq C\ep^{-2\sigma+\lambda},\notag\\
		I_\ep^2\colonequals\int\limits_{\RR_2(\GG(b_\ep,\de))}e_\ep(u,A)d\HH^2&\leq C\de^{-1}\fenthreedOmega(u,A)\leq C\ep^{-\sigma+\lambda}.\label{RR2}
	\end{align}
	Here, $C$ is a universal constant.
	
	We claim that if $\sigma<\lambda$ and $\ep$ is sufficiently small, then, for each face $\omega$ of a cube of the grid, every connected component $S_\omega$ of
	$$
	\{x\in\omega\ |\ |u(x)|\leq1/2\}
	$$
	satisfies
	$$
	d_{S_\omega}\colonequals\mathrm{deg}(u/|u|,\partial S_\omega)=0.
	$$
	Indeed, by \cite[Lemma 4.1]{roman-vortex-construction}, a result adapted from \cite{jerrard}, we have
	$$
	|d_{S_\omega}|\leq C\int_{S_\omega}|\nabla_Au|^2,
	$$
	where $C$ is a universal constant. Combining this with \eqref{RR2}, we obtain
	$$
	|d_{S_\omega}|\leq CI_\ep^2\leq C\ep^{-\sigma+\lambda}.
	$$
	Since $d_{S_\omega}$ is an integer and $\sigma<\lambda$, the claim follows for all sufficiently small $\ep$.
	
	It follows that the vorticity approximation $\nu_\ep$ of $\mu(u,A)$ provided by \cite[Theorem 1.1]{roman-vortex-construction} vanishes identically in $\overline\Omega$. Estimate \cite[(8.10)]{roman-vortex-construction} then yields
	$$
	\|\mu(u,A)\|_{C_T^{0,1}(\Omega,\R^3)^*}\leq C\de\fenthreedOmega(u,A)+C\max\{\ep\de^{-1}\fenthreedOmega(u,A),\ep\}(1+\de^{-2})\fenthreedOmega(u,A).
	$$
	Notice that
	$$
	\ep\de^{-1}\fenthreedOmega(u,A)\leq\ep\ep^{-\sigma}\ep^\lambda<\ep
	$$
	since $\sigma<\lambda$. Therefore, the previous estimate gives
	$$
	\|\mu(u,A)\|_{C_T^{0,1}(\Omega,\R^3)^*}\leq C(\de+\ep\de^{-2})\fenthreedOmega(u,A)
	$$
	and hence
	\begin{equation}\label{vorticityestimateLipschitz}
		\|\mu(u,A)\|_{C_T^{0,1}(\Omega,\R^3)^*}\leq C\ep^{\min\{\sigma,1-2\sigma\}}\fenthreedOmega(u,A).
	\end{equation}
	
	We first consider the case $\lambda>\frac13$. Choosing $\sigma=\frac13$, we have $\sigma<\lambda$ and
	$$
	\sigma\leq\frac{1+\lambda}{4}.
	$$
	Moreover,
	$$
	\min\{\sigma,1-2\sigma\}=\frac13.
	$$
	Thus, \eqref{vorticityestimateLipschitz} yields
	\begin{equation}\label{vorticityestimatelambdaabove}
		\|\mu(u,A)\|_{C_T^{0,1}(\Omega,\R^3)^*}\leq C\ep^{\frac13}\fenthreedOmega(u,A).
	\end{equation}
	
	Suppose now that $\lambda\leq\frac13$. Given $\gamma\in(0,\lambda)$, we choose
	$$
	\sigma=\lambda-\gamma.
	$$
	Then $\sigma<\lambda$, $\sigma<\frac13$, and
	$$
	\sigma\leq\frac{1+\lambda}{4}.
	$$
	Hence,
	$$
	\min\{\sigma,1-2\sigma\}=\sigma=\lambda-\gamma,
	$$
	and \eqref{vorticityestimateLipschitz} yields
	\begin{equation}\label{vorticityestimatelambdabelow}
		\|\mu(u,A)\|_{C_T^{0,1}(\Omega,\R^3)^*}\leq C\ep^{\lambda-\gamma}\fenthreedOmega(u,A).
	\end{equation}
	
	In addition, from \cite[Lemma 8.1]{roman-vortex-construction} we have
	\begin{equation}\label{vorticityestimatecontinuous}
		\|\mu(u,A)\|_{C^0(\Omega,\R^3)^*}\leq C\fenthreedOmega(u,A).
	\end{equation}
	Interpolating between \eqref{vorticityestimatelambdaabove} and \eqref{vorticityestimatecontinuous}, we conclude that, if $\lambda>\frac13$, then
	$$
	\|\mu(u,A)\|_{C_T^{0,\beta}(\Omega,\R^3)^*}\leq C_\beta\ep^{\frac{\beta}{3}}\fenthreedOmega(u,A).
	$$
	On the other hand, interpolating between \eqref{vorticityestimatelambdabelow} and \eqref{vorticityestimatecontinuous}, we conclude that, if $\lambda\leq\frac13$, then, for any $\gamma\in(0,\lambda)$,
	$$
	\|\mu(u,A)\|_{C_T^{0,\beta}(\Omega,\R^3)^*}\leq C_\beta\ep^{\beta(\lambda-\gamma)}\fenthreedOmega(u,A).
	$$
	The proposition is thus proved.
\end{proof}

%%%%%%%%%%%%%%%%%%%%%%%%%%%%%%%%%%%%%%%%%%%%%%%%%%%%%%%%%%%%%
\section{Existence}
%%%%%%%%%%%%%%%%%%%%%%%%%%%%%%%%%%%%%%%%%%%%%%%%%%%%%%%%%%%%%

We now present the proof of Theorem~\ref{theorem:meissner}.

\begin{proof}
	We proceed in several steps.
	
	\begin{enumerate}[label=\textnormal{\bfseries Step \arabic*.},ref=\arabic*,leftmargin=0pt,labelsep=*,itemindent=*,itemsep=10pt,topsep=10pt]
		
		\item {\bf Existence of a constrained minimizer.}
		
		Given a configuration $(\u,\A)$, we let $(u,A)$ be defined by the relation 
		$$
		(\u,\A)=(\meisdecom).
		$$
		We then consider the set
		$$
		U=\left\{(\u,\A)\in H^1(\Omega,\C)\times[A_\ex+H_{\curl}]\ \middle|\ \fenthreedOmega(u,A)<\ep\right\}.
		$$
		We note that $(\u_\circ,\A_\circ)\colonequals(\meisconf)$, that is, $(u_\circ,A_\circ)=(1,0)$, satisfies $\fenthreedOmega(u_\circ,A_\circ)=0$. Therefore, $U$ is non-empty.
		
		Let us prove that $\glthreed$ possesses a minimizer in $\overline U$. For this, let $\{(\tilde\u_n,\tilde\A_n)\}_n\subset U$ be a minimizing sequence. Lemma~\ref{lemma:decompositionR3} yields a gauge-transformed sequence
		$$
		\{(\u_n,\A_n)\}_n\subset H^1(\Omega,\C)\times[A_{\ex}+\dot H^1_{\diver=0}].
		$$
		By gauge invariance,
		$$
		\glthreed(\tilde\u_n,\tilde\A_n)=\glthreed(\u_n,\A_n).
		$$
		Since $(\u_\circ,\A_\circ)\in U$, we may choose the minimizing sequence so that
		$$
		\glthreed(\u_n,\A_n)\leq\glthreed(\u_\circ,\A_\circ)
		$$
		for every $n$. Then
		$$
		\|\nabla_{\A_n}\u_n\|_{L^2(\Omega,\C^3)},\quad\|a_\ep-|\u_n|^2\|_{L^2(\Omega)},\quad\mbox{and}\quad\|\curl(\A_n-A_\ex)\|_{L^2(\R^3,\R^3)}
		$$
		are bounded independently of $n$. Since $\A_n-A_\ex\in\dot H^1_{\diver=0}$, we have
		$$
		\|\A_n-A_\ex\|_{\dot H^1_{\diver=0}}=\|\curl(\A_n-A_\ex)\|_{L^2(\R^3,\R^3)},
		$$
		and therefore $\{\A_n-A_\ex\}_n$ is bounded in $\dot H^1_{\diver=0}$. Moreover, by Sobolev embedding, $\{\A_n\}_n$ is bounded in $L^6(\Omega,\R^3)$.
		
		On the other hand, using
		$$
		|\u_n|^4=(|\u_n|^2-a_\ep+a_\ep)^2\leq2(a_\ep-|\u_n|^2)^2+2a_\ep^2\leq2(a_\ep-|\u_n|^2)^2+2,
		$$
		we get
		$$
		\|\u_n\|_{L^4(\Omega,\C)}^4\leq2|\Omega|+8\ep^2\int_\Omega\frac{(a_\ep-|\u_n|^2)^2}{4\ep^2},
		$$
		which implies that $\{\u_n\}_n$ is bounded in $L^4(\Omega,\C)$. 
		
		Noting that
		$$
		\nabla\u_n=\nabla_{\A_n}\u_n+i\A_n\u_n,
		$$
		and using H\"older's inequality, we deduce that $\{\u_n\}_n$ is bounded in $H^1(\Omega,\C)$.
		
		We may then extract a subsequence, still denoted by $\{(\u_n,\A_n)\}_n$, such that
		$$
		(\u_n,\A_n-A_\ex)\rightharpoonup(\u,\A-A_\ex)
		$$
		weakly in $H^1(\Omega,\C)\times\dot H^1_{\diver=0}$. By compact Sobolev embedding, we also have strong convergence in
		$$
		L^q(\Omega,\C)\times L^q(\Omega,\R^3)
		$$
		for every $q<6$.
		
		Let us now show that
		\begin{equation}\label{lowerinfGL}
			\glthreed(\u,\A)\leq\liminf_n\glthreed(\u_n,\A_n).
		\end{equation}
		By strong convergence in $L^4(\Omega,\C)$,
		$$
		\lim_n\|a_\ep-|\u_n|^2\|_{L^2(\Omega)}=\|a_\ep-|\u|^2\|_{L^2(\Omega)}.
		$$
		Also, by the weak lower semicontinuity of the norm in $\dot H^1_{\diver=0}$,
		\begin{align*}
			\liminf_n\|\curl(\A_n-A_\ex)\|_{L^2(\R^3,\R^3)}
			&=\liminf_n\|\nabla(\A_n-A_\ex)\|_{L^2(\R^3,\R^3)}\\
			&\geq\|\nabla(\A-A_\ex)\|_{L^2(\R^3,\R^3)}\\
			&=\|\curl(\A-A_\ex)\|_{L^2(\R^3,\R^3)}.
		\end{align*}
		Moreover, $\u_n\to\u$ strongly in $L^4(\Omega,\C)$ and $\A_n\to\A$ strongly in $L^4(\Omega,\R^3)$. In particular,
		$$
		\A_n\u_n\to\A\u\quad\mbox{strongly in }L^2(\Omega,\C^3).
		$$
		Using also the weak convergence of $\nabla\u_n$ in $L^2(\Omega,\C^3)$, we deduce that
		\begin{align*}
			\liminf_n\|\nabla_{\A_n}\u_n\|_{L^2(\Omega,\C^3)}^2
			&=\liminf_n\left(\|\nabla\u_n\|_{L^2(\Omega,\C^3)}^2-2\int_\Omega(\nabla\u_n,i\A_n\u_n)+\|\A_n\u_n\|_{L^2(\Omega,\C^3)}^2\right)\\
			&\geq\|\nabla\u\|_{L^2(\Omega,\C^3)}^2-2\int_\Omega(\nabla\u,i\A\u)+\|\A\u\|_{L^2(\Omega,\C^3)}^2\\
			&=\|\nabla_\A\u\|_{L^2(\Omega,\C^3)}^2.
		\end{align*}
		Hence, \eqref{lowerinfGL} holds.
		
		We now claim that
		\begin{equation}\label{lowerinfF}
			\fenthreedOmega(u,A)\leq\liminf_n\fenthreedOmega(u_n,A_n).
		\end{equation}
		Recall that
		$$
		u_n=\pmin^{-1}e^{-ih_\ex\phi_\ep^0}\u_n,\qquad A_n=\A_n-h_\ex A_\ep^0\quad\mbox{in }\Omega.
		$$
		We first show that
		$$
		u_n\rightharpoonup u\quad\mbox{weakly in }H^1(\Omega,\C).
		$$
		Since $b\leq\pmin^2\leq1$ and $\pmin\in W^{1,\infty}(\Omega)$, multiplication by $\pmin^{-1}$ is continuous in $H^1(\Omega)$. Moreover, since $\pzero>3$, letting
		$$
		q=\frac{2\pzero}{\pzero-2}<6,
		$$
		we have
		$$
		\frac12=\frac1{\pzero}+\frac1q.
		$$
		The strong convergence of $\u_n$ to $\u$ in $L^q(\Omega,\C)$ and estimate \eqref{estimateforgradphip0} imply that
		$$
		\u_n\nabla\phi_\ep^0\longrightarrow\u\nabla\phi_\ep^0\quad\mbox{strongly in }L^2(\Omega,\C^3).
		$$
		Together with the weak convergence of $\u_n$ to $\u$ in $H^1(\Omega,\C)$, this proves the claimed weak convergence of $u_n$ in $H^1(\Omega,\C)$.
		
		On the other hand, since $A_\ep^0\in H^1(\Omega,\R^3)$ is independent of $n$, the weak convergence of $\A_n$ to $\A$ in $H^1(\Omega,\R^3)$ implies that
		$$
		A_n\rightharpoonup A\quad\mbox{weakly in }H^1(\Omega,\R^3).
		$$
		Here, we used the fact that the weak convergence of $\A_n-A_\ex$ to $\A-A_\ex$ in $\dot H^1_{\diver=0}$, together with the fact that $A_\ex$ does not depend on $n$, implies the weak convergence of $\A_n$ to $\A$ in $H^1(\Omega,\R^3)$.
		
		Arguing as in the proof of \eqref{lowerinfGL}, we deduce that
		$$
		\fenthreedOmega(u,A)\leq\liminf_n\fenthreedOmega(u_n,A_n),
		$$
		which proves \eqref{lowerinfF}.
		
		Moreover, by gauge invariance and since $(\tilde\u_n,\tilde\A_n)\in U$, we have
		$$
		\fenthreedOmega(u_n,A_n)<\ep.
		$$
		Therefore, \eqref{lowerinfF} implies that
		$$
		\fenthreedOmega(u,A)\leq\ep,
		$$
		and hence $(\u,\A)\in\overline U$.
		
		\item {\bf Interior minimality and first energy estimate.}
		
		We next show that $(\u,\A)\in U$. Since $(\u,\A)$ minimizes $\glthreed$ over $\overline U$, this will show that it is a local minimizer of $\glthreed$ and hence a critical point. Combining
		$$
		\glthreed(\u,\A)\leq\glthreed(\u_\circ,\A_\circ)=\glthreed(\meisconf)
		$$
		with \eqref{Split3D}, we find
		$$
		\fenweightthreed(u,A)-h_\ex\int_\Omega\mu(u,A)\cdot\Bfield+R_0\leq0.
		$$
		It follows that
		\begin{equation}\label{sizefen}
			\fenweightthreed(u,A)\leq h_\ex\|\mu(u,A)\|_{C_T^{0,\beta_0}(\Omega,\R^3)^*}\|\Bfield\|_{C_T^{0,\beta_0}(\Omega,\R^3)}-R_0,
		\end{equation}
		where
		\begin{equation}\label{defbeta0}
		\beta_0=1-\frac3{\pzero}>0.
		\end{equation}
		In view of \eqref{estimateforB0W1p}, Sobolev embedding gives
		\begin{equation}\label{holderbeta0}
			\|\Bfield\|_{C_T^{0,\beta_0}(\Omega,\R^3)}\leq C.
		\end{equation}
		On the other hand, since $\fenthreedOmega(u,A)\leq\ep$, Proposition~\ref{prop:vorticityestimate}, applied with $\lambda=1$, yields
		\begin{equation}\label{apriorimu}
			\|\mu(u,A)\|_{C_T^{0,\beta_0}(\Omega,\R^3)^*}\leq C\ep^{\frac{\beta_0}{3}}\fenthreedOmega(u,A).
		\end{equation}
		
		We next estimate the last term in \eqref{sizefen}. Recalling the definition of $R_0$, we have
		$$
		-R_0\leq\frac{h_\ex^2}{2b}\int_\Omega|\curl\Bfield|^2(1-|u|^2)_+.
		$$
		Let
		$$
		q=\frac{\pzero}{\pzero-2},
		$$
		so that
		$$
		\frac2{\pzero}+\frac1q=1.
		$$
		By H\"older's inequality and \eqref{estimateforB0},
		\begin{align*}
			-R_0
			&\leq Ch_\ex^2\|\curl\Bfield\|_{L^{\pzero}(\Omega,\R^3)}^2\|(1-|u|^2)_+\|_{L^q(\Omega)}\\
			&\leq Ch_\ex^2\|(1-|u|^2)_+\|_{L^q(\Omega)}.
		\end{align*}
		If $\pzero\geq4$, then $q\leq2$, and hence
		$$
		\|(1-|u|^2)_+\|_{L^q(\Omega)}
		\leq C\|(1-|u|^2)_+\|_{L^2(\Omega)}
		\leq C\ep\fenthreedOmega(u,A)^{\frac12}.
		$$
		If $3<\pzero<4$, then $q>2$. Since
		$$
		0\leq(1-|u|^2)_+\leq1,
		$$
		interpolation between $L^2(\Omega)$ and $L^\infty(\Omega)$ gives
		\begin{align*}
			\|(1-|u|^2)_+\|_{L^q(\Omega)}
			&\leq\|(1-|u|^2)_+\|_{L^2(\Omega)}^{\frac2q}\|(1-|u|^2)_+\|_{L^\infty(\Omega)}^{1-\frac2q}\\
			&\leq C\ep^{2-\frac4{\pzero}}\fenthreedOmega(u,A)^{1-\frac2{\pzero}}.
		\end{align*}
		Defining
		\begin{equation}\label{defvartheta}
		\vartheta\colonequals\min\left\{1,2-\frac4{\pzero}\right\},
		\end{equation}
		we conclude in both cases that
		\begin{equation}\label{R00}
			-R_0\leq Ch_\ex^2\ep^\vartheta\fenthreedOmega(u,A)^{\frac\vartheta2}.
		\end{equation}
		
		Combining \eqref{sizefen}, \eqref{holderbeta0}, \eqref{apriorimu}, and \eqref{R00}, and recalling that $h_\ex\leq\ep^{-\alpha}$ and $\fenthreedOmega(u,A)\leq\ep$, we obtain
		\begin{align*}
			\fenweightthreed(u,A)
			&\leq Ch_\ex\ep^{\frac{\beta_0}{3}}\fenthreedOmega(u,A)+Ch_\ex^2\ep^\vartheta\fenthreedOmega(u,A)^{\frac\vartheta2}\\
			&\leq C\ep^{1+\frac{\beta_0}{3}-\alpha}+C\ep^{\frac{3\vartheta}{2}-2\alpha}.
		\end{align*}
		Recalling \eqref{defbeta0}, we observe that the condition
		$$
		\pzero>\frac3{1-3\alpha}
		$$
		is equivalent to
		$$
		\frac{\beta_0}{3}>\alpha.
		$$
		Moreover, this condition implies that
		$$
		\frac{3\vartheta}{2}-2\alpha>1.
		$$
		Indeed, if $\pzero\geq4$, then $\vartheta=1$ and the latter inequality follows from $\alpha<\frac14$. If $3<\pzero<4$, then $\vartheta=2-\frac4{\pzero}$, and
		$$
		\frac{3\vartheta}{2}-2\alpha-1=2-\frac6{\pzero}-2\alpha>0,
		$$
		where we used $\pzero>\frac3{1-3\alpha}>\frac3{1-\alpha}$. We conclude that
		$$
		\fenweightthreed(u,A)=o(\ep).
		$$
		Since
		$$
		\fenthreedOmega(u,A)\leq b^{-2}\fenweightthreed(u,A),
		$$
		we obtain
		$$
		\fenthreedOmega(u,A)=o(\ep).
		$$
		Therefore, $(\u,\A)\in U$.
		
		Recalling \eqref{Split3D}, an immediate consequence of the previous estimates is that
		$$
		\glthreed(\u,\A)\geq\glthreed(\u_\circ,\A_\circ)+o(1).
		$$
		On the other hand, since $(\u,\A)$ minimizes $\glthreed$ over $\overline U$ and $(\u_\circ,\A_\circ)\in U$, we have
		$$
		\glthreed(\u,\A)\leq\glthreed(\u_\circ,\A_\circ).
		$$
		This proves item~\ref{part-a}.
		
		\item {\bf Clearing out and improved estimates.}
		
		Since $(\u,\A)$ is a critical point of $\glthreed$, it satisfies the Ginzburg--Landau equations. We can gauge transform this configuration into $(\uCouGauge,\ACouGauge)$ in the Coulomb gauge in $\Omega$. Notice that
		$$
		(\pmin^{-1}\uCouGauge,\ACouGauge)
		$$
		is gauge equivalent to
		$$
		(u,A+h_\ex(A_\ep^0-\nabla\phi_\ep^0)).
		$$
		Hence, using \eqref{defBfield}, $\|\curl\Bfield\|_{L^2(\Omega,\R^3)}\leq C$, and $\|u\|_{L^\infty(\Omega,\C)}\leq1$, we find
		\begin{align*}
			\|\nabla_{\ACouGauge}(\pmin^{-1}\uCouGauge)\|^2_{L^2(\Omega,\C^3)}
			&=\|\nabla_{A+h_\ex(A_\ep^0-\nabla\phi_\ep^0)}u\|^2_{L^2(\Omega,\C^3)}\\
			&=\|\nabla_{A+h_\ex\pmin^{-2}\curl\Bfield}u\|^2_{L^2(\Omega,\C^3)}\\
			&\leq C\left(\|\nabla_Au\|^2_{L^2(\Omega,\C^3)}+h_\ex^2\right)\\
			&\leq C(\ep+h_\ex^2).
		\end{align*}
		Proposition~\ref{prop:regularityA} then implies that
		$$
		\|\ACouGauge\|_{L^\infty(\Omega,\R^3)}\leq C\ep^{-\alpha}.
		$$
		Using \eqref{boundnablau} and the fact that $\alpha<1$, we find
		$$
		\|\nabla\uCouGauge\|_{L^\infty(\Omega,\C^3)}\leq C\ep^{-1}.
		$$
		By gauge invariance, $|\u|=|\uCouGauge|$, and therefore
		$$
		\|\nabla|\u|\|_{L^\infty(\Omega,\R^3)}=\|\nabla|\uCouGauge|\|_{L^\infty(\Omega,\R^3)}\leq\|\nabla\uCouGauge\|_{L^\infty(\Omega,\C^3)}\leq\frac{C}{\ep}.
		$$
		Moreover, since
		$$
		|u|=\frac{|\u|}{\pmin},
		$$
		Proposition~\ref{prop:regularitypmin} yields
		$$
		\|\nabla|u|\|_{L^\infty(\Omega,\R^3)}=\left\|\nabla\frac{|\u|}{\pmin}\right\|_{L^\infty(\Omega,\R^3)}\leq C\left(\|\nabla|\u|\|_{L^\infty(\Omega,\R^3)}+\|\nabla\pmin\|_{L^\infty(\Omega,\R^3)}\right)\leq\frac{C}{\ep}.
		$$
		We can therefore apply Lemma~\ref{lemma:clearingout}. For every $x_0\in\overline\Omega$,
		$$
		(1-|u|(x_0))^5\leq\frac{C}{\ep^3}\int_{B_\ep(x_0)\cap\Omega}(1-|u|^2)^2\leq\frac{C}{\ep^3}\int_\Omega(1-|u|^2)^2\leq\frac{C}{\ep}\fenthreedOmega(u,A)=o(1).
		$$
		We have thus proved that
		$$
		\|1-|u|\|_{L^\infty(\Omega)}=o(1).
		$$
		Since $|\u|=\pmin|u|$, this also implies that
		$$
		\||\u|-\pmin\|_{L^\infty(\Omega)}=o(1).
		$$
		This proves item~\ref{part-b}.
		
		We can now use \cite[Proposition A.1]{Rom-CMP} to improve our vorticity estimate. Since $|u|\geq\frac12$ for all sufficiently small $\ep$, we have
		\begin{equation}\label{vortfinal}
			\|\mu(u,A)\|_{C_T^{0,1}(\Omega,\R^3)^*}\leq C\ep\fenthreedOmega(u,A).
		\end{equation}
		By interpolation with \eqref{vorticityestimatecontinuous}, we find, for any $\beta\in(0,1]$,
		$$
		\|\mu(u,A)\|_{C_T^{0,\beta}(\Omega,\R^3)^*}\leq C\ep^\beta\fenthreedOmega(u,A).
		$$
		This allows us to improve the estimate on $\fenthreedOmega(u,A)$. Arguing as above, with $\beta=\beta_0$, and recalling \eqref{defbeta0} and \eqref{defvartheta}, we obtain
		\begin{equation}\label{improvedenergy}
			\fenweightthreed(u,A)\leq Ch_\ex\ep^{\beta_0}\fenthreedOmega(u,A)+Ch_\ex^2\ep^\vartheta\fenthreedOmega(u,A)^{\frac\vartheta2}.
		\end{equation}
		Since $\beta_0>3\alpha>\alpha$, the first term on the right-hand side can be absorbed for all sufficiently small $\ep$. Since $\frac\vartheta2\leq\frac12$, using once again
		$$
		\fenthreedOmega(u,A)\leq b^{-2}\fenweightthreed(u,A),
		$$
		we deduce from \eqref{improvedenergy} that
		$$
		\fenthreedOmega(u,A)^{1-\frac\vartheta2}\leq C\ep^{\vartheta-2\alpha}.
		$$
		Therefore,
		\begin{equation*}
			\fenthreedOmega(u,A)\leq C\ep^{\frac{\vartheta-2\alpha}{1-\frac\vartheta2}}.
		\end{equation*}
		Let
		$$
		\de\colonequals\frac{\vartheta-2\alpha}{1-\frac\vartheta2}-1=\frac{\frac{3\vartheta}{2}-2\alpha-1}{1-\frac\vartheta2}.
		$$
		The inequality
		$$
		\frac{3\vartheta}{2}-2\alpha>1
		$$
		proved above shows that $\de>0$. Hence,
		$$
		\fenthreedOmega(u,A)\leq C\ep^{1+\de}.
		$$
		This proves item~\ref{part-c}.
		
		On the other hand, by direct computation,
		\begin{align*}
			\mu\left(\frac{\u}{\pmin},\A\right)
			&=\mu\left(\frac{\u}{\pmin},\A-h_\ex(A_\ep^0-\nabla\phi_\ep^0)+h_\ex(A_\ep^0-\nabla\phi_\ep^0)\right)\\
			&=\mu\left(\frac{\u}{\pmin},\A-h_\ex(A_\ep^0-\nabla\phi_\ep^0)\right)+h_\ex\curl\left(\left(1-\left|\frac{\u}{\pmin}\right|^2\right)(A_\ep^0-\nabla\phi_\ep^0)\right)\\
			&=\mu(u,A)+h_\ex\curl\left((1-|u|^2)\frac{\curl\Bfield}{\pmin^2}\right).
		\end{align*}
		Here, in the last equality, we used that
		$$
		\left(\pmin^{-1}\u,\A-h_\ex(A_\ep^0-\nabla\phi_\ep^0)\right)
		$$
		is gauge equivalent to $(u,A)$, together with \eqref{defBfield}.
		
		\smallskip
		For any $B\in C_T^{0,1}(\Omega,\R^3)$, integration by parts gives
		\begin{align*}
			\left|\int_\Omega h_\ex\curl\left((1-|u|^2)\frac{\curl\Bfield}{\pmin^2}\right)\cdot B\right|
			&=\left|h_\ex\int_\Omega(1-|u|^2)\frac{\curl\Bfield}{\pmin^2}\cdot\curl B\right|\\
			&\leq Ch_\ex\|1-|u|^2\|_{L^2(\Omega)}\|\curl\Bfield\|_{L^2(\Omega,\R^3)}\|B\|_{C_T^{0,1}(\Omega,\R^3)}\\
			&\leq Ch_\ex\ep\fenthreedOmega(u,A)^{\frac12}\|B\|_{C_T^{0,1}(\Omega,\R^3)}\\
			&=o(1)\|B\|_{C_T^{0,1}(\Omega,\R^3)}.
		\end{align*}
		We therefore deduce that
		$$
		\left\|\mu\left(\frac{\u}{\pmin},\A\right)-\mu(u,A)\right\|_{C_T^{0,1}(\Omega,\R^3)^*}=o(1).
		$$
		Together with \eqref{vortfinal}, this proves item~\ref{part-d}.
		
		\item {\bf Convergence in the global divergence-free gauge.}
		
		By gauge invariance and Lemma~\ref{lemma:decompositionR3}, we may assume that
		$$
		\A-A_\ex\in\dot H^1_{\diver=0}.
		$$
		Recalling that
		$$
		A=\A-h_\ex A_\ep^0=(\A-A_\ex)-h_\ex(A_\ep^0-A_{0,\ex}),
		$$
		and using \eqref{spaceforA0}, we deduce that
		$$
		A\in\dot H^1_{\diver=0}.
		$$
		Moreover,
		\begin{equation}\label{smallnessA}
		\|A\|_{\dot H^1_{\diver=0}}^2=\|\curl A\|_{L^2(\R^3,\R^3)}^2\leq2\fenweightthreed(u,A)=o(\ep).
		\end{equation}
		On the other hand,
		\begin{align*}
			\int_\Omega|\nabla u|^2
			&\leq2\int_\Omega|\nabla_Au|^2+2\int_\Omega|u|^2|A|^2\\
			&\leq C\int_\Omega\pmin^2|\nabla_Au|^2+C\|u\|_{L^\infty(\Omega,\C)}^2\|A\|_{L^2(\Omega,\R^3)}^2\\
			&\leq C\fenweightthreed(u,A)+C\|A\|_{\dot H^1_{\diver=0}}^2\\
			&\leq C\fenweightthreed(u,A)=o(\ep),
		\end{align*}
		where we used $b\leq\pmin^2$, $\|1-|u|\|_{L^\infty(\Omega)}=o(1)$, and
		$$
		\|A\|_{L^2(\Omega,\R^3)}\leq C\|A\|_{L^6(\R^3,\R^3)}\leq C\|A\|_{\dot H^1_{\diver=0}}.
		$$
		Letting
		$$
		\underline u=\frac1{|\Omega|}\int_\Omega u,
		$$
		the Poincar\'e--Wirtinger inequality yields
		\begin{equation}\label{derivate}
			\int_\Omega|u-\underline u|^2\leq C\int_\Omega|\nabla u|^2=o(\ep).
		\end{equation}
		Notice also that
		$$
		\int_\Omega(1-|u|)^2\leq\int_\Omega(1-|u|^2)^2
		$$
		and
		$$
		\int_\Omega\left||u|-|\underline u|\right|^2\leq\int_\Omega|u-\underline u|^2.
		$$
		We then deduce that
		$$
		\int_\Omega|1-|\underline u||^2\leq C\int_\Omega(1-|u|^2)^2+C\int_\Omega|u-\underline u|^2=o(\ep).
		$$
		Since $\underline u$ is constant, it follows that
		$$
		\underline u=(1+o(\sqrt\ep))e^{i\theta}=e^{i\theta}+o(\sqrt\ep)
		$$
		for some $\theta=\theta(\ep)\in[0,2\pi]$. From \eqref{derivate}, we obtain
		\begin{equation}\label{combine0}
			\inf_{\theta\in[0,2\pi]}\|u-e^{i\theta}\|_{H^1(\Omega,\C)}=o(\sqrt\ep).
		\end{equation}
		Combining this with \eqref{smallnessA}, we have proved item~\ref{part-e}.
		
		\item {\bf Convergence to the Meissner configuration.}
		
		Recalling that
		$$
		\pmin^{-1}\u=u\phase,
		$$
		we have
		\begin{equation}\label{combine1}
			\inf_{\theta\in[0,2\pi]}\|\pmin^{-1}\u-e^{i\theta}\phase\|_{H^1(\Omega,\C)}=\inf_{\theta\in[0,2\pi]}\|(u-e^{i\theta})\phase\|_{H^1(\Omega,\C)}.
		\end{equation}
		On the other hand, for any $\theta\in[0,2\pi]$,
		\begin{equation}\label{combine2}
			\|(u-e^{i\theta})\phase\|_{H^1(\Omega,\C)}\leq C\left(\|u-e^{i\theta}\|_{L^2(\Omega,\C)}+\|\nabla u\|_{L^2(\Omega,\C^3)}+h_\ex\|(u-e^{i\theta})\nabla\phi_\ep^0\|_{L^2(\Omega,\C^3)}\right).
		\end{equation}
		Recalling \eqref{estimateforgradphip0}, since $\pzero>3$, H\"older's inequality and Sobolev embedding yield
		\begin{equation}\label{combine3}
			\|(u-e^{i\theta})\nabla\phi_\ep^0\|_{L^2(\Omega,\C^3)}
			\leq\|u-e^{i\theta}\|_{L^{\frac{2\pzero}{\pzero-2}}(\Omega,\C)}\|\nabla\phi_\ep^0\|_{L^{\pzero}(\Omega,\R^3)}
			\leq C\|u-e^{i\theta}\|_{H^1(\Omega,\C)}.
		\end{equation}
		Since $h_\ex\leq\ep^{-\alpha}$ and $\alpha<\frac14$, \eqref{combine0} and \eqref{combine3} imply that
		$$
		h_\ex\inf_{\theta\in[0,2\pi]}\|(u-e^{i\theta})\nabla\phi_\ep^0\|_{L^2(\Omega,\C^3)}=o(1).
		$$
		Combining this with \eqref{combine1} and \eqref{combine2}, we finally obtain
		$$
		\inf_{\theta\in[0,2\pi]}\|\pmin^{-1}\u-e^{i\theta}\phase\|_{H^1(\Omega,\C)}=o(1).
		$$
		Since
		$$
		\A-h_\ex A_\ep^0=A,
		$$
		and $\|A\|_{\dot H^1_{\diver=0}}=o(1)$, item~\ref{part-f} follows.
		
	\end{enumerate}
	
	The theorem is thus proved.
\end{proof}

\begin{remark}
	The condition $\pzero>\frac3{1-3\alpha}$ arises from the preliminary estimate of the vorticity term. Indeed, Proposition~\ref{prop:vorticityestimate} provides a gain of order $\ep^{\frac{\beta}{3}}$ in the dual of $C_T^{0,\beta}(\Omega,\R^3)$, while
	$$
	W^{1,\pzero}(\Omega)\hookrightarrow C^{0,1-\frac3{\pzero}}(\overline\Omega).
	$$
	Taking $\beta=1-\frac3{\pzero}$ gives the condition
	$$
	\frac13\left(1-\frac3{\pzero}\right)>\alpha.
	$$
	The natural gain in the vorticity estimate would be of order $\ep^\beta$. This would yield the condition
	$$
	\pzero>\frac3{1-\alpha},
	$$
	which coincides with the condition appearing in the uniqueness result and seems to be the natural, and possibly sharp, threshold. Obtaining such a gain, however, appears to be very difficult within the available three-dimensional $\ep$-level vortex construction.
\end{remark}

%%%%%%%%%%%%%%%%%%%%%%%%%%%%%%%%%%%%%%%%%%%%%%%%%%%%%%%%%%%%%
\section{Uniqueness}
%%%%%%%%%%%%%%%%%%%%%%%%%%%%%%%%%%%%%%%%%%%%%%%%%%%%%%%%%%%%%

We now present the proof of Theorem~\ref{theorem:uniqueness}.

\begin{proof}
	We proceed in several steps.
	
	\begin{enumerate}[label=\textnormal{\bfseries Step \arabic*.},ref=\arabic*,leftmargin=0pt,labelsep=*,itemindent=*,itemsep=10pt,topsep=10pt]
		
		\item {\bf Coulomb gauge estimates and clearing out.}
		
		Let us assume towards a contradiction that there exist two local minimizers $(\u_j,\A_j)$ satisfying the hypotheses of the theorem which are not gauge equivalent. Writing
		$$
		(\u_j,\A_j)=(\pmin u_j\phase,A_j+h_\ex A_\ep^0),
		$$
		we see that this configuration is gauge equivalent to
		$$
		(\pmin u_j,A_j+h_\ex(A_\ep^0-\nabla\phi_\ep^0)).
		$$
		We gauge transform the latter configuration into $(\uCouGauge_j,\ACouGauge_j)$, where $\ACouGauge_j$ is in the Coulomb gauge in $\Omega$. Thus,
		\begin{equation*}
			\left\{
			\begin{array}{rcll}
				\diver\ACouGauge_j&=&0&\mathrm{in}\ \Omega\\
				\ACouGauge_j\cdot\nu&=&0&\mathrm{on}\ \partial\Omega.
			\end{array}
			\right.
		\end{equation*}
		In particular,
		$$
		|\uCouGauge_j|=\pmin|u_j|.
		$$
		
		Using the gauge invariance of $\glthreed$, we deduce that $(\uCouGauge_j,\ACouGauge_j)$ is a critical point. Proposition~\ref{prop:regularity} then yields
		$$
		\|\uCouGauge_j\|_{L^\infty(\Omega,\C)}\leq1\quad\mbox{and}\quad\left\|\frac{\uCouGauge_j}{\pmin}\right\|_{L^\infty(\Omega,\C)}\leq1.
		$$
		
		Since $(\uCouGauge_j,\ACouGauge_j)$ is in the Coulomb gauge in $\Omega$, Proposition~\ref{prop:regularityA} gives
		$$
		\|\ACouGauge_j\|_{L^\infty(\Omega,\R^3)}\leq C\left(\left\|\nabla_{\ACouGauge_j}\left(\frac{\uCouGauge_j}{\pmin}\right)\right\|_{L^2(\Omega,\C^3)}+h_\ex\right).
		$$
		Moreover,
		$$
		\left(\frac{\uCouGauge_j}{\pmin},\ACouGauge_j\right)
		$$
		is gauge equivalent to
		$$
		(u_j,A_j+h_\ex(A_\ep^0-\nabla\phi_\ep^0)).
		$$
		Therefore, by the gauge invariance of the covariant gradient,
		\begin{align*}
			\left\|\nabla_{\ACouGauge_j}\left(\frac{\uCouGauge_j}{\pmin}\right)\right\|_{L^2(\Omega,\C^3)}^2
			&=\|\nabla_{A_j+h_\ex(A_\ep^0-\nabla\phi_\ep^0)}u_j\|_{L^2(\Omega,\C^3)}^2\\
			&\leq C\left(\|\nabla_{A_j}u_j\|_{L^2(\Omega,\C^3)}^2+h_\ex^2\|u_j(A_\ep^0-\nabla\phi_\ep^0)\|_{L^2(\Omega,\C^3)}^2\right)\\
			&\leq C\left(\fenthreedOmega(u_j,A_j)+h_\ex^2\right).
		\end{align*}
		Here, we used \eqref{defBfield}, together with $\|u_j\|_{L^\infty(\Omega,\C)}\leq1$ and $\|\curl\Bfield\|_{L^2(\Omega,\R^3)}\leq C$. Since
		$$
		\fenthreedOmega(u_j,A_j)\leq c_0\ep
		$$
		and $h_\ex\leq\ep^{-\alpha}$, we conclude that
		\begin{equation}\label{boundAWLinfty}
			\|\ACouGauge_j\|_{L^\infty(\Omega,\R^3)}\leq C\ep^{-\alpha}.
		\end{equation}
		
		In turn, Proposition~\ref{prop:regularity} and the fact that $\alpha<1$ yield
		$$
		\|\nabla\uCouGauge_j\|_{L^\infty(\Omega,\C^3)}\leq C\left(\ep^{-1}+\|\ACouGauge_j\|_{L^\infty(\Omega,\R^3)}\right)\leq C\ep^{-1}.
		$$
		Using
		$$
		|\nabla|\uCouGauge_j||\leq|\nabla\uCouGauge_j|
		$$
		and
		$$
		|u_j|=\pmin^{-1}|\uCouGauge_j|,
		$$
		together with Proposition~\ref{prop:regularitypmin}, we find
		\begin{equation*}
			\|\nabla|u_j|\|_{L^\infty(\Omega,\R^3)}=\|\nabla(\pmin^{-1}|\uCouGauge_j|)\|_{L^\infty(\Omega,\R^3)}\leq C\ep^{-1}.
		\end{equation*}
		
		We can now apply Lemma~\ref{lemma:clearingout}. Denoting by $C_{\mathrm{co}}>0$ the constant in the clearing-out estimate, for every $x_0\in\overline\Omega$ we have
		$$
		(1-|u_j|(x_0))^5\leq\frac{C_{\mathrm{co}}}{\ep^3}\int_{B_\ep(x_0)\cap\Omega}(1-|u_j|^2)^2\leq\frac{C_{\mathrm{co}}}{\ep}\fenthreedOmega(u_j,A_j)\leq C_{\mathrm{co}}c_0.
		$$
		Let
		\begin{equation*}
			\tau\colonequals(C_{\mathrm{co}}c_0)^{\frac15}.
		\end{equation*}
		We choose $c_0>0$ sufficiently small so that
		\begin{equation}\label{choicec0}
			\tau<\frac14\quad\mbox{and}\quad4\tau\leq\frac{\de_0}{2},
		\end{equation}
		where $\de_0=\de_0(b,\pzero,r_0)$ is the constant in \eqref{conditionBMO}. In particular, $c_0$ depends on the fixed parameters of the problem and on $\pzero$.
		
		By the definition of $\tau$, we have
		$$
		\frac34\leq1-\tau\leq|u_j|\leq1\quad\mathrm{in}\ \Omega.
		$$
		We may thus write
		$$
		\frac{\uCouGauge_j}{\pmin}=\eta_j e^{i\phi_j},
		$$
		where
		\begin{equation}\label{boundeta}
			\eta_j=\left|\frac{\uCouGauge_j}{\pmin}\right|=|u_j|,
			\quad\mbox{with }\frac34\leq1-\tau\leq\eta_j\leq1.
		\end{equation}
		
		\item {\bf Elimination of the phase.}
		
		We gauge transform $(\uCouGauge_j,\ACouGauge_j)$ into $(\uUniGauge_j,\AUniGauge_j)$, where
		$$
		\uUniGauge_j=\pmin\eta_j,\qquad\AUniGauge_j=\ACouGauge_j-\nabla\phi_j.
		$$
		After extending $\phi_j$ to a function in $H^1_{\loc}(\R^3)$, this gauge transformation is defined in the whole space. In particular, $(\uUniGauge_j,\AUniGauge_j)$ is a critical point of $\glthreed$ and is gauge equivalent to $(\u_j,\A_j)$.
		
		Since $(\uUniGauge_j,\AUniGauge_j)$ satisfies the second equation in \eqref{GLeq}, we have
		$$
		\curl^2(\AUniGauge_j-A_\ex)=(i\uUniGauge_j,\nabla_{\AUniGauge_j}\uUniGauge_j)\chi_\Omega\quad\mathrm{in}\ \R^3.
		$$
		Taking the divergence of this equation in $\Omega$, we find
		$$
		\diver\big((i\uUniGauge_j,\nabla_{\AUniGauge_j}\uUniGauge_j)\big)=0\quad\mathrm{in}\ \Omega.
		$$
		Since $\uUniGauge_j=\pmin\eta_j$ is real-valued,
		$$
		(i\uUniGauge_j,\nabla_{\AUniGauge_j}\uUniGauge_j)=-\pmin^2\eta_j^2\AUniGauge_j\quad\mathrm{in}\ \Omega.
		$$
		Hence,
		$$
		\diver(\eta_j^2\pmin^2\AUniGauge_j)=0\quad\mathrm{in}\ \Omega.
		$$
		The first boundary condition in \eqref{GLeq} also yields
		$$
		\AUniGauge_j\cdot\nu=0\quad\mathrm{on}\ \partial\Omega,
		$$
		because $\uUniGauge_j$ is real-valued and strictly positive. Recalling that
		$$
		\AUniGauge_j=\ACouGauge_j-\nabla\phi_j,
		$$
		we conclude that
		\begin{equation}\label{pdephi}
			\left\{
			\begin{array}{rcll}
				\diver\big(\eta_j^2\pmin^2(\ACouGauge_j-\nabla\phi_j)\big)&=&0&\mathrm{in}\ \Omega\\
				(\ACouGauge_j-\nabla\phi_j)\cdot\nu&=&0&\mathrm{on}\ \partial\Omega.
			\end{array}
			\right.
		\end{equation}
		
		Thus, $\ACouGauge_j-\nabla\phi_j$ is in the weighted Coulomb gauge associated with the weight $\eta_j^2\pmin^2$. This is the only point in the paper where we need to apply the weighted elliptic estimate to a perturbation of the original weight $\pmin^2$. The factor $\frac12$ in \eqref{conditionBMO} was introduced precisely to ensure that the small-BMO condition required in \cite[Theorem 2]{SunWang} remains valid after this perturbation.
		
		Indeed, by assumption \eqref{conditionBMO},
		$$
		[\pmin^2]_{BMO_{r_0}(\Omega)}\leq\frac{\de_0}{2}.
		$$
		Moreover, \eqref{boundeta} implies that
		$$
		\|\eta_j^2-1\|_{L^\infty(\Omega)}\leq2\tau.
		$$
		For any $f\in L^\infty(\Omega)$, we have
		$$
		[f]_{BMO_{r_0}(\Omega)}\leq2\|f\|_{L^\infty(\Omega)}.
		$$
		Therefore, using $\pmin^2\leq1$, we obtain
		\begin{equation*}
			[(\eta_j^2-1)\pmin^2]_{BMO_{r_0}(\Omega)}
			\leq2\|(\eta_j^2-1)\pmin^2\|_{L^\infty(\Omega)}
			\leq2\|\eta_j^2-1\|_{L^\infty(\Omega)}.
		\end{equation*}
		Consequently,
		\begin{align*}
			[\eta_j^2\pmin^2]_{BMO_{r_0}(\Omega)}\leq[\pmin^2]_{BMO_{r_0}(\Omega)}+[(\eta_j^2-1)\pmin^2]_{BMO_{r_0}(\Omega)}			&\leq\frac{\de_0}{2}+2\|\eta_j^2-1\|_{L^\infty(\Omega)}\\
			&\leq\frac{\de_0}{2}+4\tau \leq \de_0,
		\end{align*}
		where in the last inequality we used \eqref{choicec0}.
		
		Therefore, \cite[Theorem 2]{SunWang} applied to \eqref{pdephi} yields
		$$
		\|\nabla\phi_j\|_{L^{\pzero}(\Omega,\R^3)}\leq C\|\ACouGauge_j\|_{L^{\pzero}(\Omega,\R^3)}.
		$$
		Combining this with \eqref{boundAWLinfty}, we find
		$$
		\|\nabla\phi_j\|_{L^{\pzero}(\Omega,\R^3)}\leq C\ep^{-\alpha}.
		$$
		It follows that
		\begin{equation}\label{boundAUni}
			\|\AUniGauge_j\|_{L^{\pzero}(\Omega,\R^3)}\leq\|\ACouGauge_j\|_{L^{\pzero}(\Omega,\R^3)}+\|\nabla\phi_j\|_{L^{\pzero}(\Omega,\R^3)}\leq C\ep^{-\alpha}.
		\end{equation}
		
		\item {\bf Strict energy comparison.}
		
		We prove that, for all sufficiently small $\ep$,
		\begin{equation}\label{energycomparison}
			\frac{Y}{2}\colonequals\frac{\glthreed(\uUniGauge_1,\AUniGauge_1)+\glthreed(\uUniGauge_2,\AUniGauge_2)}2-\glthreed\left(\frac{\uUniGauge_1+\uUniGauge_2}2,\frac{\AUniGauge_1+\AUniGauge_2}2\right)>0.
		\end{equation}
		Recalling that $\uUniGauge_j=\pmin\eta_j$ and using \eqref{lm mag decoupling}, we find
		$$
		\glthreed(\uUniGauge_j,\AUniGauge_j)=E_\ep(\pmin)+\frac12\int_\Omega\left(\pmin^2|\nabla\eta_j|^2+\pmin^2\eta_j^2|\AUniGauge_j|^2+\frac{\pmin^4}{\ep^2}(1-\eta_j^2)^2\right)+\frac12\int_{\R^3}|\curl\AUniGauge_j-H_\ex|^2.
		$$
		We write
		$$
		Y=Y_0+Y_1+Y_2+Y_3,
		$$
		where
		\begin{align*}
			Y_0&=\frac12\int_\Omega\pmin^2\left(|\nabla\eta_1|^2+|\nabla\eta_2|^2\right)-\int_\Omega\pmin^2\left|\nabla\left(\frac{\eta_1+\eta_2}2\right)\right|^2\\
			Y_1&=\frac12\int_\Omega\pmin^2\left(\eta_1^2|\AUniGauge_1|^2+\eta_2^2|\AUniGauge_2|^2\right)-\int_\Omega\pmin^2\left(\frac{\eta_1+\eta_2}2\right)^2\left|\frac{\AUniGauge_1+\AUniGauge_2}2\right|^2\\
			Y_2&=\frac1{4\ep^2}\int_\Omega\pmin^4\left((1-\eta_1^2)^2+(1-\eta_2^2)^2\right)-\frac1{2\ep^2}\int_\Omega\pmin^4\left(1-\left(\frac{\eta_1+\eta_2}2\right)^2\right)^2\\
			Y_3&=\frac12\int_{\R^3}\left(|\curl\AUniGauge_1-H_\ex|^2+|\curl\AUniGauge_2-H_\ex|^2\right)-\int_{\R^3}\left|\curl\left(\frac{\AUniGauge_1+\AUniGauge_2}2\right)-H_\ex\right|^2.
		\end{align*}
		Observe that $Y_3\geq 0$ by convexity. We next prove that $Y_0+Y_1+Y_2>0$. A direct computation gives
		\begin{equation}\label{Y0}
			Y_0=\frac14\int_\Omega\pmin^2|\nabla(\eta_1-\eta_2)|^2\geq\frac b4\|\nabla(\eta_1-\eta_2)\|_{L^2(\Omega,\R^3)}^2.
		\end{equation}
		Moreover,
		\begin{multline}\label{Y1}
			Y_1=\int_\Omega\pmin^2\Bigg(\frac18(\eta_1-\eta_2)^2\left(|\AUniGauge_1|^2+|\AUniGauge_2|^2\right)+\frac1{16}(3\eta_1^2-2\eta_1\eta_2+3\eta_2^2)|\AUniGauge_1-\AUniGauge_2|^2\\
			+\frac18(\eta_1-\eta_2)(\AUniGauge_1-\AUniGauge_2)\cdot\left((\eta_1+3\eta_2)\AUniGauge_1+(3\eta_1+\eta_2)\AUniGauge_2\right)\Bigg).
		\end{multline}
		Finally,
		$$
		Y_2=\frac1{32\ep^2}\int_\Omega\pmin^4(\eta_1-\eta_2)^2\left(7(\eta_1+\eta_2)^2-4\eta_1\eta_2-8\right).
		$$
		In view of \eqref{boundeta}, we have
		$$
		7(\eta_1+\eta_2)^2-4\eta_1\eta_2-8\geq\frac{11}{2}.
		$$
		Using also $\pmin^4\geq b^2$, we conclude that
		\begin{equation}\label{Y2}
			Y_2\geq\frac{11b^2}{64\ep^2}\|\eta_1-\eta_2\|_{L^2(\Omega)}^2.
		\end{equation}
		
		We now estimate the last term in \eqref{Y1}. Since $\pmin^2\leq1$ and $\eta_j\leq1$, its absolute value is bounded by
		$$
		C\int_\Omega|\eta_1-\eta_2||\AUniGauge_1-\AUniGauge_2|\left(|\AUniGauge_1|+|\AUniGauge_2|\right).
		$$
		Let
		$$
		r=\frac{2\pzero}{\pzero-2},
		$$
		so that
		$$
		\frac1r+\frac12+\frac1{\pzero}=1.
		$$
		Therefore, by H\"older's inequality and \eqref{boundAUni},
		\begin{align*}
			&\int_\Omega|\eta_1-\eta_2||\AUniGauge_1-\AUniGauge_2|\left(|\AUniGauge_1|+|\AUniGauge_2|\right)\\
			&\qquad\leq\|\eta_1-\eta_2\|_{L^r(\Omega)}\|\AUniGauge_1-\AUniGauge_2\|_{L^2(\Omega,\R^3)}\left(\|\AUniGauge_1\|_{L^{\pzero}(\Omega,\R^3)}+\|\AUniGauge_2\|_{L^{\pzero}(\Omega,\R^3)}\right)\\
			&\qquad\leq C\ep^{-\alpha}\|\eta_1-\eta_2\|_{L^r(\Omega)}\|\AUniGauge_1-\AUniGauge_2\|_{L^2(\Omega,\R^3)}.
		\end{align*}
		Since $\pzero>3$, we have $2<r<6$. Interpolation and Sobolev embedding give
		\begin{align*}
			\|\eta_1-\eta_2\|_{L^r(\Omega)}
			&\leq C\|\eta_1-\eta_2\|_{L^2(\Omega)}^{1-\frac3{\pzero}}\|\eta_1-\eta_2\|_{L^6(\Omega)}^{\frac3{\pzero}}\\
			&\leq C\|\eta_1-\eta_2\|_{L^2(\Omega)}^{1-\frac3{\pzero}}\left(\|\nabla(\eta_1-\eta_2)\|_{L^2(\Omega,\R^3)}+\|\eta_1-\eta_2\|_{L^2(\Omega)}\right)^{\frac3{\pzero}}.
		\end{align*}
		Writing
		$$
		\|\eta_1-\eta_2\|_{L^2(\Omega)}=\ep\left(\ep^{-2}\|\eta_1-\eta_2\|_{L^2(\Omega)}^2\right)^{\frac12}
		$$
		and observing that, for $\ep<1$,
		\begin{multline*}
			\|\nabla(\eta_1-\eta_2)\|_{L^2(\Omega,\R^3)}+\|\eta_1-\eta_2\|_{L^2(\Omega)}\\
			\leq C\left(\|\nabla(\eta_1-\eta_2)\|_{L^2(\Omega,\R^3)}^2+\ep^{-2}\|\eta_1-\eta_2\|_{L^2(\Omega)}^2\right)^{\frac12},
		\end{multline*}
		it follows that
		$$
		\|\eta_1-\eta_2\|_{L^r(\Omega)}
		\leq C\ep^{1-\frac3{\pzero}}\left(\|\nabla(\eta_1-\eta_2)\|_{L^2(\Omega,\R^3)}^2+\ep^{-2}\|\eta_1-\eta_2\|_{L^2(\Omega)}^2\right)^{\frac12}.
		$$
		Consequently,
		\begin{multline*}
			\int_\Omega|\eta_1-\eta_2||\AUniGauge_1-\AUniGauge_2|\left(|\AUniGauge_1|+|\AUniGauge_2|\right)\\
			\leq C\ep^{1-\alpha-\frac3{\pzero}}\left(\|\nabla(\eta_1-\eta_2)\|_{L^2(\Omega,\R^3)}^2+\ep^{-2}\|\eta_1-\eta_2\|_{L^2(\Omega)}^2\right)^{\frac12}\|\AUniGauge_1-\AUniGauge_2\|_{L^2(\Omega,\R^3)}.
		\end{multline*}
		By Young's inequality,
		\begin{multline}\label{mixedterm}
			\int_\Omega|\eta_1-\eta_2||\AUniGauge_1-\AUniGauge_2|\left(|\AUniGauge_1|+|\AUniGauge_2|\right)\\
			\leq C\ep^{1-\alpha-\frac3{\pzero}}\left(\|\nabla(\eta_1-\eta_2)\|_{L^2(\Omega,\R^3)}^2+\ep^{-2}\|\eta_1-\eta_2\|_{L^2(\Omega)}^2+\|\AUniGauge_1-\AUniGauge_2\|_{L^2(\Omega,\R^3)}^2\right).
		\end{multline}
		Since
		$$
		\pzero>\frac3{1-\alpha},
		$$
		we have
		$$
		\ep^{1-\alpha-\frac3{\pzero}}=o(1).
		$$
		
		In view of \eqref{boundeta}, we have
		$$
		3\eta_1^2-2\eta_1\eta_2+3\eta_2^2=(\eta_1-\eta_2)^2+2\eta_1^2+2\eta_2^2\geq\frac94.
		$$
		Therefore, the second term in \eqref{Y1} satisfies
		\begin{equation}\label{positiveG}
			\frac1{16}\int_\Omega\pmin^2(3\eta_1^2-2\eta_1\eta_2+3\eta_2^2)|\AUniGauge_1-\AUniGauge_2|^2\geq\frac{9b}{64}\|\AUniGauge_1-\AUniGauge_2\|_{L^2(\Omega,\R^3)}^2.
		\end{equation}
		Combining \eqref{Y0}, \eqref{Y1}, \eqref{Y2}, \eqref{mixedterm}, and \eqref{positiveG}, and discarding the first nonnegative term in \eqref{Y1}, we obtain
		\begin{multline*}
			Y_0+Y_1+Y_2\geq\left(\frac b4-o(1)\right)\|\nabla(\eta_1-\eta_2)\|_{L^2(\Omega,\R^3)}^2\\
			+\left(\frac{11b^2}{64}-o(1)\right)\ep^{-2}\|\eta_1-\eta_2\|_{L^2(\Omega)}^2+\left(\frac{9b}{64}-o(1)\right)\|\AUniGauge_1-\AUniGauge_2\|_{L^2(\Omega,\R^3)}^2.
		\end{multline*}
		Therefore,
		$$
		Y_0+Y_1+Y_2>0
		$$
		for all sufficiently small $\ep$, unless
		$$
		\eta_1=\eta_2\quad\mbox{and}\quad\AUniGauge_1=\AUniGauge_2.
		$$
		Since the two configurations are not gauge equivalent, the latter alternative is impossible. Recalling that $Y_3\geq0$, we conclude that
		$$
		Y>0.
		$$
		
		\item {\bf Contradiction.}
		
		We may assume without loss of generality that
		$$
		\glthreed(\uUniGauge_1,\AUniGauge_1)\leq\glthreed(\uUniGauge_2,\AUniGauge_2).
		$$
		Let
		$$
		(\u^{(1)},\A^{(1)})\colonequals\left(\frac{\uUniGauge_1+\uUniGauge_2}{2},\frac{\AUniGauge_1+\AUniGauge_2}{2}\right).
		$$
		From \eqref{energycomparison}, we have
		$$
		\glthreed(\u^{(1)},\A^{(1)})<\frac{\glthreed(\uUniGauge_1,\AUniGauge_1)+\glthreed(\uUniGauge_2,\AUniGauge_2)}2\leq\glthreed(\uUniGauge_2,\AUniGauge_2).
		$$
		
		We define inductively
		$$
		(\u^{(n+1)},\A^{(n+1)})\colonequals\left(\frac{\u^{(n)}+\uUniGauge_2}{2},\frac{\A^{(n)}+\AUniGauge_2}{2}\right).
		$$
		Each configuration $(\u^{(n)},\A^{(n)})$ is of the form
		$$
		(\u^{(n)},\A^{(n)})=(\pmin\eta^{(n)},\A^{(n)}),
		$$
		where, by construction and in view of \eqref{boundeta},
		$$
		\frac34\leq1-\tau\leq\eta^{(n)}\leq1.
		$$
		Moreover, \eqref{boundAUni} yields
		$$
		\|\A^{(n)}\|_{L^{\pzero}(\Omega,\R^3)}\leq C\ep^{-\alpha}
		$$
		uniformly in $n$. The argument of the previous step therefore applies to the pair $(\u^{(n)},\A^{(n)})$ and $(\uUniGauge_2,\AUniGauge_2)$ and gives
		$$
		\glthreed(\u^{(n+1)},\A^{(n+1)})<\frac{\glthreed(\u^{(n)},\A^{(n)})+\glthreed(\uUniGauge_2,\AUniGauge_2)}2<\glthreed(\uUniGauge_2,\AUniGauge_2).
		$$
		It follows by induction that
		$$
		\glthreed(\u^{(n)},\A^{(n)})<\glthreed(\uUniGauge_2,\AUniGauge_2)
		$$
		for every $n\geq1$.
		
		On the other hand,
		$$
		(\u^{(n)},\A^{(n)})-(\uUniGauge_2,\AUniGauge_2)=2^{-n}\left((\uUniGauge_1,\AUniGauge_1)-(\uUniGauge_2,\AUniGauge_2)\right).
		$$
		Therefore,
		$$
		(\u^{(n)},\A^{(n)})\longrightarrow(\uUniGauge_2,\AUniGauge_2)
		$$
		with respect to the natural topology of $H^1(\Omega,\C)\times[A_\ex+H_{\curl}]$. This contradicts the local minimality of $(\uUniGauge_2,\AUniGauge_2)$. Hence,
		$$
		(\uUniGauge_1,\AUniGauge_1)=(\uUniGauge_2,\AUniGauge_2).
		$$
		
	\end{enumerate}
	
	The theorem is thus proved.
\end{proof}

\bibliography{references3DPinning}

% \bib, bibdiv, biblist are defined by the amsrefs package.
\begin{bibdiv}
\begin{biblist}

\bib{amrouche-seloula}{article}{
      author={Amrouche, Ch\'erif},
      author={Seloula, Nour El~Houda},
       title={{$L^p$}-theory for vector potentials and {S}obolev's inequalities
  for vector fields: application to the {S}tokes equations with pressure
  boundary conditions},
        date={2013},
        ISSN={0218-2025,1793-6314},
     journal={Math. Models Methods Appl. Sci.},
      volume={23},
      number={1},
       pages={37\ndash 92},
         url={https://doi.org/10.1142/S0218202512500455},
      review={\MR{2997467}},
}

\bib{BBH-CVPDE}{article}{
      author={Bethuel, Fabrice},
      author={Brezis, Ha\"im},
      author={H\'elein, Fr\'ed\'eric},
       title={Asymptotics for the minimization of a {G}inzburg-{L}andau
  functional},
        date={1993},
        ISSN={0944-2669,1432-0835},
     journal={Calc. Var. Partial Differential Equations},
      volume={1},
      number={2},
       pages={123\ndash 148},
         url={https://doi.org/10.1007/BF01191614},
      review={\MR{1261720}},
}

\bib{BetBreHel}{book}{
      author={Bethuel, Fabrice},
      author={Brezis, Ha\"im},
      author={H\'elein, Fr\'ed\'eric},
       title={Ginzburg-{L}andau vortices},
      series={Progress in Nonlinear Differential Equations and their
  Applications},
   publisher={Birkh\"auser Boston, Inc., Boston, MA},
        date={1994},
      volume={13},
        ISBN={0-8176-3723-0},
         url={http://dx.doi.org/10.1007/978-1-4612-0287-5},
      review={\MR{1269538}},
}

\bib{BetBreOrl}{article}{
      author={Bethuel, F.},
      author={Brezis, H.},
      author={Orlandi, G.},
       title={Asymptotics for the {G}inzburg-{L}andau equation in arbitrary
  dimensions},
        date={2001},
        ISSN={0022-1236,1096-0783},
     journal={J. Funct. Anal.},
      volume={186},
      number={2},
       pages={432\ndash 520},
         url={https://doi.org/10.1006/jfan.2001.3791},
      review={\MR{1864830}},
}

\bib{BCS}{article}{
      author={Bardeen, J.},
      author={Cooper, L.~N.},
      author={Schrieffer, J.~R.},
       title={Theory of superconductivity},
        date={1957},
     journal={Phys. Rev.},
      volume={108},
      number={5},
       pages={1175\ndash 1204},
         url={https://link.aps.org/doi/10.1103/PhysRev.108.1175},
}

\bib{BalJerOrlSon1}{article}{
      author={Baldo, S.},
      author={Jerrard, R.~L.},
      author={Orlandi, G.},
      author={Soner, H.~M.},
       title={Convergence of {G}inzburg-{L}andau functionals in
  three-dimensional superconductivity},
        date={2012},
        ISSN={0003-9527},
     journal={Arch. Ration. Mech. Anal.},
      volume={205},
      number={3},
       pages={699\ndash 752},
         url={http://dx.doi.org/10.1007/s00205-012-0527-2},
      review={\MR{2960031}},
}

\bib{BalJerOrlSon2}{article}{
      author={Baldo, S.},
      author={Jerrard, R.~L.},
      author={Orlandi, G.},
      author={Soner, H.~M.},
       title={Vortex density models for superconductivity and superfluidity},
        date={2013},
        ISSN={0010-3616},
     journal={Comm. Math. Phys.},
      volume={318},
      number={1},
       pages={131\ndash 171},
         url={http://dx.doi.org/10.1007/s00220-012-1629-2},
      review={\MR{3017066}},
}

\bib{SunWang}{article}{
      author={Byun, Sun-Sig},
      author={Wang, Lihe},
       title={The conormal derivative problem for elliptic equations with {BMO}
  coefficients on {R}eifenberg flat domains},
        date={2005},
        ISSN={0024-6115,1460-244X},
     journal={Proc. London Math. Soc. (3)},
      volume={90},
      number={1},
       pages={245\ndash 272},
         url={https://doi.org/10.1112/S0024611504014960},
      review={\MR{2107043}},
}

\bib{chapman}{article}{
      author={Chapman, S.~Jonathan},
       title={Superheating field of type {II} superconductors},
        date={1995},
        ISSN={0036-1399},
     journal={SIAM J. Appl. Math.},
      volume={55},
      number={5},
       pages={1233\ndash 1258},
         url={https://doi.org/10.1137/S0036139993254760},
      review={\MR{1349308}},
}

\bib{phys-x-2024}{article}{
      author={Chen, C.},
      author={Liu, Y.},
      author={Chen, Y.},
      author={Hu, Y.~N.},
      author={Zhang, T.~Z.},
      author={Li, D.},
      author={Wang, X.},
      author={Wang, C.~X.},
      author={Lu, Z. Y.~W.},
      author={Zhang, Y.~H.},
      author={Zhang, Q.~L.},
      author={Dong, X.~L.},
      author={Wang, R.},
      author={Feng, D.~L.},
      author={Zhang, T.},
       title={Revealing the microscopic mechanism of elementary vortex pinning
  in superconductors},
        date={2024},
     journal={Phys. Rev. X},
      volume={14},
       pages={041039},
         url={https://link.aps.org/doi/10.1103/PhysRevX.14.041039},
}

\bib{chapman-richardson-pinning}{article}{
      author={Chapman, S.~J.},
      author={Richardson, G.},
       title={Vortex pinning by inhomogeneities in type-{II} superconductors},
        date={1997},
        ISSN={0167-2789,1872-8022},
     journal={Phys. D},
      volume={108},
      number={4},
       pages={397\ndash 407},
         url={https://doi.org/10.1016/S0167-2789(97)00053-5},
      review={\MR{1474691}},
}

\bib{DosSantos3D}{article}{
      author={Dos~Santos, Micka\"el},
       title={Study of a {$3D$}-{G}inzburg-{L}andau functional with a
  discontinuous pinning term},
        date={2012},
        ISSN={0362-546X,1873-5215},
     journal={Nonlinear Anal.},
      volume={75},
      number={17},
       pages={6275\ndash 6296},
         url={https://doi.org/10.1016/j.na.2012.07.004},
      review={\MR{2959805}},
}

\bib{DosRodSan}{article}{
      author={Dos~Santos, Micka\"el},
      author={Rodiac, R\'emy},
      author={Sandier, Etienne},
       title={The {G}inzburg-{L}andau energy with a pinning term oscillating
  faster than the coherence length},
        date={2023},
        ISSN={1463-9963,1463-9971},
     journal={Interfaces Free Bound.},
      volume={25},
      number={3},
       pages={491\ndash 515},
         url={https://doi.org/10.4171/ifb/495},
      review={\MR{4642021}},
}

\bib{diaz-roman-3d}{article}{
      author={D\'iaz-Vera, Mat\'ias},
      author={Rom\'an, Carlos},
       title={First critical field in the pinned three-dimensional
  {G}inzburg-{L}andau model of superconductivity},
        date={2026},
        ISSN={0944-2669,1432-0835},
     journal={Calc. Var. Partial Differential Equations},
      volume={65},
      number={7},
       pages={Paper No. 217, 37},
         url={https://doi.org/10.1007/s00526-026-03383-5},
      review={\MR{5090628}},
}

\bib{diaz-roman-2d}{article}{
      author={D\'iaz-Vera, Mat\'ias},
      author={Rom\'an, Carlos},
       title={On the {M}eissner state for type-{II} inhomogeneous
  superconductors},
        date={2026},
        ISSN={0022-1236,1096-0783},
     journal={J. Funct. Anal.},
      volume={290},
      number={7},
       pages={Paper No. 111332, 45},
         url={https://doi.org/10.1016/j.jfa.2025.111332},
      review={\MR{5013119}},
}

\bib{DRY}{unpublished}{
      author={D{\'i}az-Vera, Mat{\'i}as},
      author={Rom{\'a}n, Carlos},
      author={Yamamoto, Kai},
        note={Work in progress},
}

\bib{GinLan}{article}{
      author={Ginzburg, V.~L.},
      author={Landau, L.~D.},
       title={On the theory of superconductivity},
        date={1950},
     journal={Zh. Eksp. Teor. Fiz.},
      volume={20},
       pages={1064\ndash 1082},
        note={English translation in: \emph{Collected papers of L.D.Landau},
  Edited by D. Ter. Haar, Pergamon Press, Oxford 1965, pp. 546--568},
}

\bib{Groger}{article}{
      author={Gr\"oger, Konrad},
       title={A {$W^{1,p}$}-estimate for solutions to mixed boundary value
  problems for second order elliptic differential equations},
        date={1989},
        ISSN={0025-5831,1432-1807},
     journal={Math. Ann.},
      volume={283},
      number={4},
       pages={679\ndash 687},
         url={https://doi.org/10.1007/BF01442860},
      review={\MR{990595}},
}

\bib{jerrard}{article}{
      author={Jerrard, Robert~L.},
       title={Lower bounds for generalized {G}inzburg-{L}andau functionals},
        date={1999},
        ISSN={0036-1410,1095-7154},
     journal={SIAM J. Math. Anal.},
      volume={30},
      number={4},
       pages={721\ndash 746},
         url={https://doi.org/10.1137/S0036141097300581},
      review={\MR{1684723}},
}

\bib{KozSoh}{article}{
      author={Kozono, Hideo},
      author={Sohr, Hermann},
       title={New a priori estimates for the {S}tokes equations in exterior
  domains},
        date={1991},
        ISSN={0022-2518},
     journal={Indiana Univ. Math. J.},
      volume={40},
      number={1},
       pages={1\ndash 27},
         url={http://dx.doi.org/10.1512/iumj.1991.40.40001},
      review={\MR{1101219}},
}

\bib{Kwok2016}{article}{
      author={Kwok, Wai-Kwong},
      author={Welp, Ulrich},
      author={Glatz, Andreas},
      author={Koshelev, Alexei~E.},
      author={Kihlstrom, Karen~J.},
      author={Crabtree, George~W.},
       title={Vortices in high-performance high-temperature superconductors},
        date={2016},
        ISSN={0034-4885,1361-6633},
     journal={Rep. Prog. Phys.},
      volume={79},
      number={11},
       pages={116501, 39},
         url={https://doi.org/10.1088/0034-4885/79/11/116501},
}

\bib{Likha}{article}{
      author={Likharev, K.~K.},
       title={Superconducting weak links},
        date={1979},
     journal={Reviews of Modern Physics},
      volume={51},
      number={1},
       pages={101\ndash 159},
}

\bib{lassoued-mironescu}{article}{
      author={Lassoued, Lotfi},
      author={Mironescu, Petru},
       title={Ginzburg-{L}andau type energy with discontinuous constraint},
        date={1999},
        ISSN={0021-7670,1565-8538},
     journal={J. Anal. Math.},
      volume={77},
       pages={1\ndash 26},
         url={https://doi.org/10.1007/BF02791255},
      review={\MR{1753480}},
}

\bib{LinRiviere}{article}{
      author={Lin, Fanghua},
      author={Rivi\`ere, Tristan},
       title={Complex {G}inzburg-{L}andau equations in high dimensions and
  codimension two area minimizing currents},
        date={1999},
        ISSN={1435-9855,1435-9863},
     journal={J. Eur. Math. Soc. (JEMS)},
      volume={1},
      number={3},
       pages={237\ndash 311},
         url={https://doi.org/10.1007/s100970050008},
      review={\MR{1714735}},
}

\bib{Meyers}{article}{
      author={Meyers, Norman~G.},
       title={An {$L\sp{p}$}-estimate for the gradient of solutions of second
  order elliptic divergence equations},
        date={1963},
        ISSN={0391-173X},
     journal={Ann. Scuola Norm. Sup. Pisa Cl. Sci. (3)},
      volume={17},
       pages={189\ndash 206},
      review={\MR{159110}},
}

\bib{Riv}{article}{
      author={Rivi\`ere, Tristan},
       title={Line vortices in the {${\rm U}(1)$}-{H}iggs model},
        date={1995/96},
        ISSN={1292-8119,1262-3377},
     journal={ESAIM Contr\^ole Optim. Calc. Var.},
      volume={1},
       pages={77\ndash 167},
         url={https://doi.org/10.1051/cocv:1996103},
      review={\MR{1394302}},
}

\bib{Rom-CMP}{article}{
      author={Rom\'an, Carlos},
       title={On the first critical field in the three dimensional
  {G}inzburg-{L}andau model of superconductivity},
        date={2019},
        ISSN={0010-3616,1432-0916},
     journal={Comm. Math. Phys.},
      volume={367},
      number={1},
       pages={317\ndash 349},
         url={https://doi.org/10.1007/s00220-019-03306-w},
      review={\MR{3933412}},
}

\bib{roman-vortex-construction}{article}{
      author={Rom\'an, Carlos},
       title={Three dimensional vortex approximation construction and
  {$\varepsilon$}-level estimates for the {G}inzburg-{L}andau functional},
        date={2019},
        ISSN={0003-9527,1432-0673},
     journal={Arch. Ration. Mech. Anal.},
      volume={231},
      number={3},
       pages={1531\ndash 1614},
         url={https://doi.org/10.1007/s00205-018-1304-7},
      review={\MR{3902469}},
}

\bib{RomanPinningUpper25}{article}{
      author={Rom\'an, Carlos},
       title={First critical field in the pinned three-dimensional
  {Ginzburg--Landau} model: A matching upper bound},
        date={2026},
     journal={Ann. Henri Poincaré},
        note={Available at \url{https://doi.org/10.1007/s00023-026-01754-4}},
}

\bib{roman-sandier-serfaty}{article}{
      author={Rom\'an, Carlos},
      author={Sandier, Etienne},
      author={Serfaty, Sylvia},
       title={Bounded vorticity for the 3{D} {G}inzburg-{L}andau model and an
  isoflux problem},
        date={2023},
        ISSN={0024-6115,1460-244X},
     journal={Proc. Lond. Math. Soc. (3)},
      volume={126},
      number={3},
       pages={1015\ndash 1062},
         url={https://doi.org/10.1112/plms.12505},
      review={\MR{4563866}},
}

\bib{RSS2}{unpublished}{
      author={Rom\'{a}n, Carlos},
      author={Sandier, Etienne},
      author={Serfaty, Sylvia},
       title={Vortex lines interaction in the three-dimensional magnetic
  {G}inzburg--{Landau} model},
        date={2025},
        note={Available at
  \href{https://arxiv.org/abs/2510.14910}{arXiv:2510.14910}},
}

\bib{sandier-surface-ball-construction}{article}{
      author={Sandier, Etienne},
       title={Ginzburg-{L}andau minimizers from {$\mathbb{R}^{n+1}$} to
  {$\mathbb{R}^n$} and minimal connections},
        date={2001},
        ISSN={0022-2518,1943-5258},
     journal={Indiana Univ. Math. J.},
      volume={50},
      number={4},
       pages={1807\ndash 1844},
         url={https://doi.org/10.1512/iumj.2001.50.1751},
      review={\MR{1889083}},
}

\bib{Ser2}{article}{
      author={Serfaty, Sylvia},
       title={Stable configurations in superconductivity: uniqueness,
  multiplicity, and vortex-nucleation},
        date={1999},
        ISSN={0003-9527,1432-0673},
     journal={Arch. Ration. Mech. Anal.},
      volume={149},
      number={4},
       pages={329\ndash 365},
         url={https://doi.org/10.1007/s002050050177},
      review={\MR{1731999}},
}

\bib{SanSerBook}{book}{
      author={Sandier, Etienne},
      author={Serfaty, Sylvia},
       title={Vortices in the magnetic {G}inzburg-{L}andau model},
      series={Progress in Nonlinear Differential Equations and their
  Applications},
   publisher={Birkh\"auser Boston, Inc., Boston, MA},
        date={2007},
      volume={70},
        ISBN={978-0-8176-4316-4; 0-8176-4316-8},
      review={\MR{2279839}},
}

\bib{Tin}{book}{
      author={Tinkham, Michael},
       title={Introduction to superconductivity},
     edition={Second},
   publisher={McGraw-Hill, New York},
        date={1996},
}

\end{biblist}
\end{bibdiv}
\end{document}